\documentclass[preprint]{elsarticle}

\usepackage{jcomp}
\usepackage{framed,multirow}

\usepackage{amssymb}
\usepackage{latexsym}
\usepackage{slashbox}

\usepackage{url}
\usepackage{color}
\usepackage{lineno}
\usepackage{amsmath,mathrsfs}
\usepackage[caption=false]{subfig}
\usepackage{dsfont}
\usepackage{float}
\usepackage{gensymb}
\usepackage{natbib}
\usepackage{graphicx,titlesec}
\usepackage{epstopdf}
\usepackage{geometry}
\usepackage{booktabs}
\definecolor{newcolor}{rgb}{.8,.349,.1}

\usepackage{varioref}
\usepackage{hyperref}
\usepackage{cleveref}

\hypersetup{
	colorlinks=true,
	linkcolor={red!50!black},
	citecolor={blue!50!black},
	urlcolor={blue!80!black}
}

\newcommand\dt {{\Delta t}}
\def\m{\mbox{\boldmath $m$}}
\def\bM{\mbox{\boldmath $M$}}
\def\bx{\mbox{\boldmath $x$}}
\def\by{\mbox{\boldmath $y$}}

\newtheorem{example}{Example}[section]

\crefformat{equation}{(#2#1#3)}
\crefmultiformat{equation}{(#2#1#3)}{ and~(#2#1#3)}{, (#2#1#3)}{ and~(#2#1#3)}

\crefformat{figure}{Figure~#2#1#3}
\crefmultiformat{figure}{Figures~(#2#1#3)}{ and~(#2#1#3)}{, (#2#1#3)}{ and~(#2#1#3)}
\crefformat{table}{Table~#2#1#3}

\begin{document}


\begin{frontmatter}
\title{Two improved Gauss-Seidel projection methods for Landau-Lifshitz-Gilbert equation}

\author[1]{Panchi Li}
\ead{LiPanchi1994@163.com}
\author[1]{Changjian Xie}

\ead{20184007005@stu.suda.edu.cn}
\author[1,2]{Rui Du\corref{cor1}}
\cortext[cor1]{Corresponding authors}
\ead{durui@suda.edu.cn}
\author[1,2]{Jingrun Chen\corref{cor1}}
\ead{jingrunchen@suda.edu.cn}
\author[3]{Xiao-Ping Wang\corref{cor1}}
\ead{mawang@ust.hk}

\address[1]{School of Mathematical Sciences, Soochow University, Suzhou, 215006, China.}
\address[2]{Mathematical Center for Interdisciplinary Research, Soochow University, Suzhou, 215006, China.}
\address[3]{Department of Mathematics, The Hong Kong University of Science and Technology, Clear Water Bay, Kowloon, Hong Kong, China}

\begin{abstract}
Micromagnetic simulation is an important tool to study various dynamic behaviors of magnetic order in ferromagnetic materials.
The underlying model is the Landau-Lifshitz-Gilbert equation, where the magnetization dynamics is driven by
the gyromagnetic torque term and the Gilbert damping term. Numerically, considerable progress has been made in the past decades.
One of the most popular methods is the Gauss-Seidel projection method developed by Xiao-Ping Wang, Carlos Garc\'{i}a-Cervera, and Weinan E in 2001.
It first solves a set of heat equations with constant coefficients and updates the gyromagnetic term in the Gauss-Seidel manner,
and then solves another set of heat equations with constant coefficients for the damping term.
Afterwards, a projection step is applied to preserve the length constraint in the pointwise sense.
This method has been verified to be unconditionally stable numerically and successfully applied to study magnetization dynamics
under various controls.

In this paper, we present two improved Gauss-Seidel projection methods with unconditional stability.
The first method updates the gyromagnetic term and the damping term simultaneously and follows
by a projection step. The second method introduces two sets of approximate solutions, where
we update the gyromagnetic term and the damping term simultaneously for one set of approximate
solutions and apply the projection step to the other set of approximate solutions in an alternating manner.
Compared to the original Gauss-Seidel projection method which has to solve heat equations $7$ times at each time step,
the improved methods solve heat equations $5$ times and $3$ times, respectively.
First-order accuracy in time and second-order accuracy in space are verified by examples in both 1D and 3D.
In addition, unconditional stability with respect to both the grid size and the damping parameter is confirmed
numerically. Application of both methods to a realistic material is also presented with hysteresis loops
and magnetization profiles. Compared with the original method, the recorded running times suggest that savings
of both methods are about $2/7$ and $4/7$ for the same accuracy requirement, respectively.
\end{abstract}

\begin{keyword}

\KWD Landau-Lifshitz-Gilbert equation\sep Gauss-Seidel projection method\sep unconditional stability\sep micromagnetic simulation\\
\MSC[2000] 35Q99 \sep 65Z05 \sep 65M06
\end{keyword}

\end{frontmatter}


\section{Introduction}

In ferromagnetic materials, the intrinsic magnetic order, known as magnetization $\bM = (M_1, M_2, M_3)^T$,  is modeled by
the following Landau-Lifshitz-Gilbert (LLG) equation \cite{LandauLifshitz1935, Gilbert1955, WBrown1963}
\begin{equation}
\frac {\partial\mathbf{M}}{\partial t} = -\gamma\mathbf{M}\times\boldsymbol{\mathcal{H}} -
\frac {\gamma\alpha}{M_s}\mathbf{M}\times(\mathbf{M}\times\boldsymbol{\mathcal{H}})
\label{LL-Equation}
\end{equation}
with $\gamma$ the gyromagnetic ratio and $|\mathbf{M}| = M_s$ the saturation magnetization.
On the right-hand side of \eqref{LL-Equation}, the first term is the gyromagnetic
term and the second term is the Gilbert damping term with $\alpha$ the dimensionless
damping coefficient \cite{Gilbert1955}. Note that the gyromagnetic term is a conservative
term, whereas the damping term is a dissipative term.
The local field $\boldsymbol{\mathcal{H}} = -\frac{\delta F}{\delta\mathbf{M}}$ is computed from the Landau-Lifshitz energy functional
\begin{equation}
\label{LL-Energy}
F[\mathbf{M}] = \frac 1{2}\int_\Omega\left\{ \frac{A}{M_s^2}|\nabla\mathbf{M}|^2 +
\Phi\left(\frac{\mathbf{M}}{M_s}\right)
 -2\mu_0\mathbf{H}_e\cdot\mathbf{M} \right\}\mathrm{d}\bx
+ \frac{\mu_0}{2}\int_{\mathbb{R}^3}|\nabla U|^2 \mathrm{d}\bx,
\end{equation}
where $A$ is the exchange constant, $\frac{A}{M_s^2}|\nabla\mathbf{M}|$ is the
exchange interaction energy; $\Phi\left(\frac{\mathbf{M}}{M_s}\right)$
is the anisotropy energy, and for simplicity the material is assumed to be uniaxial with
$\Phi\left(\frac{\mathbf{M}}{M_s}\right) = \frac{K_u}{M_s^2}(M_2^2+M_3^2)$ with
$K_u$ the anisotropy constant; $-2\mu_0\mathbf{H}_e\cdot\mathbf{M}$ is the Zeeman energy due to the
external field with $\mu_0$ the permeability of vacuum. $\Omega$ is the volume occupied by the material.
The last term in \eqref{LL-Energy} is the energy resulting from the field
induced by the magnetization distribution inside the material.
This stray field $\boldsymbol{\mathbf{H}}_s = -\nabla U$ where $U(\bx)$ satisfies
\begin{equation}
U(\bx) = \int_{\Omega}\nabla N(\bx-\by)\cdot\mathbf{M}(\by)\mathrm{d}\by,
\end{equation}
where $N(\bx-\by) = -\frac 1{4\pi}\frac 1{|\bx-\by|}$ is the Newtonian potential.

For convenience, we rescale the original LLG equation \eqref{LL-Equation} by changes of variables
$t \rightarrow(\mu_0\gamma M_s)^{-1}t$ and  $x\rightarrow Lx$ with $L$ the diameter of $\Omega$.
Define $\mathbf{m}=\bM/M_s$ and $\mathbf{h} = M_s\mathcal{H}$. The dimensionless LLG equation reads as
\begin{equation}
\label{LLDimensionless}
\frac{\partial\mathbf{m}}{\partial t} = -\mathbf{m}\times\mathbf{h} -
\alpha\mathbf{m}\times(\mathbf{m}\times\mathbf{h}),
\end{equation}
where
\begin{equation}
\mathbf{h} = -Q(m_2\mathbf{e}_2 + m_3\mathbf{e}_3) + \epsilon\Delta\mathbf{m}
+ \mathbf{h}_e + \mathbf{h}_s
\label{dimensionlessH}
\end{equation}
with dimensionless parameters $Q=K_u/(\mu_0M_s^2)$ and $\epsilon=A/(\mu_0M_s^2L^2)$.
Here $\mathbf{e}_2=(0,1,0)$, $\mathbf{e}_3=(0,0,1)$.
Neumann boundary condition is used
\begin{equation}
\frac{\partial\mathbf{m}}{\partial\boldsymbol{\nu}}_{|\partial\Omega} = 0,
\label{boundary}
\end{equation}
where $\boldsymbol{\nu}$ is the outward unit normal vector on $\partial\Omega$.

The LLG equation is a weakly nonlinear equation. In the absence of Gilbert damping,
$\alpha = 0$, equation \eqref{LLDimensionless} is a degenerate equation of parabolic type and is
related to the sympletic flow of harmonic maps \cite{Sulem1986}. In the large
damping limit, $\alpha \to \infty$, equation \eqref{LLDimensionless} is related to
the heat flow for harmonic maps \cite{Struwe1988}. It is easy to check that $\lvert \mathbf{m}\rvert = 1$
in the pointwise sense in the evolution. All these properties possesses interesting
challenges for designing numerical methods to solve the LLG equation.
Meanwhile, micromagnetic simulation is an important tool to study magnetization
dynamics of magnetic materials \cite{WBrown1963, ZuticFabianDasSarma:2004}.
Over the past decades, there has been increasing progress on numerical methods
for the LLG equation; see \cite{Prohl2006,Cimrak2007,CJreview2007} for reviews and
references therein. Finite difference method and finite element method have been
used for the spatial discretization.

For the temporal discretization, there are explicit schemes such as Runge-Kutta
methods \cite{alouges2006convergence,FourRK2008}. Their stepsizes are subject to strong stability constraint.
Another issue is that the length of magnetization cannot be preserved and thus a projection
step is needed. Implicit schemes \cite{Yamada2004Implicit, bartels2006convergence, implicit2012}
are unconditionally stable and usually can preserve the length of magnetization automatically.
The difficulty of implicit schemes is how to solve a nonlinear system of equations at each step.
Therefore, semi-implicit methods \cite{NumGSPM2001,NumMethods2000,SecSemi2019,ImproveGSPM2003,seim2005}
provide a compromise between stability and the difficult for solving the equation at each step.
A projection step is also needed to preserve the length of magnetization.

Among the semi-implicit schemes, the most popular one is the Gauss-Seidel projection method (GSPM) proposed
by Wang, Garc\'{i}a-Cervera, and E \cite{NumGSPM2001,ImproveGSPM2003}.
GSPM first solves a set of heat equations with constant coefficients and updates the gyromagnetic term in the Gauss-Seidel manner,
and then solves another set of heat equations with constant coefficients for the damping term.
Afterwards, a projection step is applied to preserve the length of magnetization.
GSPM is first-order accurate in time and has been verified to be unconditionally stable numerically.

In this paper, we present two improved Gauss-Seidel projection methods with unconditional stability.
The first method updates the gyromagnetic term and the damping term simultaneously and follows
by a projection step. The second method introduces two sets of approximate solutions, where
we update the gyromagnetic term and the damping term simultaneously for one set of approximate
solutions and apply the projection step to the other set of approximate solutions in an alternating manner.
Compared to the original Gauss-Seidel projection method, which solves heat equations $7$ times at each time step,
the improved methods solve heat equations $5$ times and $3$ times, respectively.
First-order accuracy in time and second-order accuracy in space are verified by examples in both 1D and 3D.
In addition, unconditional stability with respect to both the grid size and the damping parameter is confirmed
numerically. Application of both methods to a realistic material is also presented with hysteresis loops
and magnetization profiles. Compared with the original method, the recorded running times suggest that savings
of both methods are about $2/7$ and $4/7$ for the same accuracy requirement, respectively.

The rest of the paper is organized as follows. For completeness and comparison, we first introduce
GSPM in \Cref{SectionGSPM}. Two improved GSPMs are presented in \Cref{SectionImproved}.
Detailed numerical tests are given in \Cref{sectionNumerical}, including accuracy check and efficiency check
in both 1D and 3D, unconditional stability with respect to both the grid size and the damping parameter,
hysteresis loops, and magnetization profiles. Conclusions are drawn in \Cref{sectionClusion}.

\section{Gauss-Seidel projection method for Landau-Lifshitz-Gilbert equation}\label{SectionGSPM}

\begin{figure}[ht]
	\centering
	\includegraphics[width = 14cm]{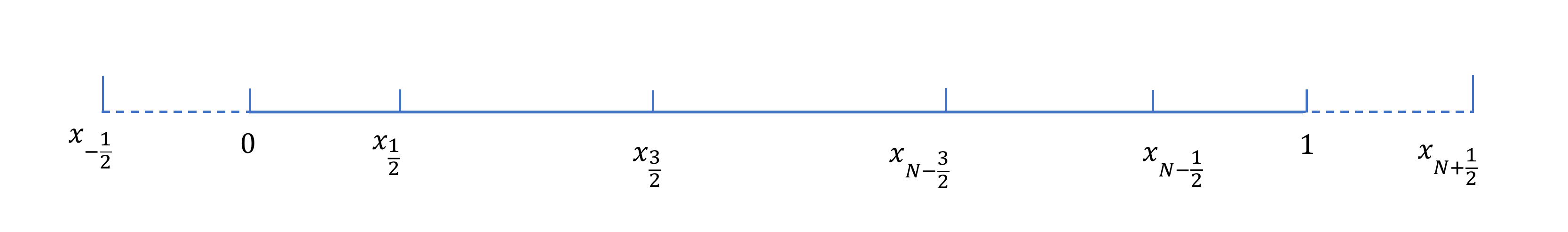}
	\caption{Spatial grids in 1D. Nodes $x_{-\frac1{2}}$ and $x_{N+\frac1{2}}$ are ghost points.}
	\label{meshgrid}
\end{figure}

Before the introduction of the GSPM \cite{NumGSPM2001,ImproveGSPM2003},
we first use the finite difference method for spatial discretization. Figure \ref{meshgrid} shows
a schematic picture of spatial grids in 1D. Let $i=0,1,\cdots,M,M+1$,
$j=0,1,\cdots,N,N+1$, and $k=0,1,\cdots,K,K+1$ be the indices of grid points in 3D.

Second-order centered difference
for $\Delta\mathbf{m}$ reads as
\begin{align}
\Delta_h\mathbf{m}_{i,j,k} &=\frac{\mathbf{m}_{i+1,j,k}-
	2\mathbf{m}_{i,j,k}+\mathbf{m}_{i-1,j,k}}{\Delta x^2}\nonumber \\
&+\frac{\mathbf{m}_{i,j+1,k}-2\mathbf{m}_{i,j,k}+
	\mathbf{m}_{i,j-1,k}}{\Delta y^2}\nonumber \\
&+\frac{\mathbf{m}_{i,j,k+1}-2\mathbf{m}_{i,j,k}+
	\mathbf{m}_{i,j,k-1}}{\Delta z^2},
\label{LaplacianDiscrete}
\end{align}
where $\mathbf{m}_{i,j,k}=\mathbf{m}((i - \frac 1{2})\Delta x,
(j - \frac 1{2})\Delta y, (k - \frac 1{2})\Delta z)$. For the Neumann boundary condition, a second-order approximation yields
\begin{align*}
\mathbf{m}_{0,j,k} & =\mathbf{m}_{1,j,k},\quad \mathbf{m}_{M,j,k}  =\mathbf{m}_{M+1,j,k},\quad j = 1,\cdots,N,k=1,\cdots,K, \\
\mathbf{m}_{i,0,k} & =\mathbf{m}_{i,1,k},\quad \mathbf{m}_{i,N,k}  =\mathbf{m}_{i,N+1,k},\quad i = 1,\cdots,M,k=1,\cdots,K, \\
\mathbf{m}_{i,j,0} & =\mathbf{m}_{i,j,1},\quad \mathbf{m}_{i,j,K}  =\mathbf{m}_{i,j,K+1},\quad i = 1,\cdots,M,j=1,\cdots,N.
\end{align*}

To illustrate the main ideas, we first consider the following simplified equation
\begin{equation*}
\mathbf{m}_t = - \mathbf{m}\times \Delta\mathbf{m} -
\alpha\mathbf{m}\times(\mathbf{m}\times\Delta\mathbf{m}),
\end{equation*}
which can be rewritten as
\begin{equation}
\label{AFLL}
\mathbf{m}_t = - \mathbf{m}\times \Delta\mathbf{m} - \alpha\mathbf{m}
(\mathbf{m}\cdot\Delta\mathbf{m}) + \alpha\Delta\mathbf{m}.
\end{equation}

We split \eqref{AFLL} into two equations
\begin{align}
\mathbf{m}_t & = - \mathbf{m}\times \Delta\mathbf{m}, \label{eqn:splitting1}\\
\mathbf{m}_t & =  \alpha\Delta\mathbf{m}.\label{eqn:splitting2}
\end{align}

However, \eqref{eqn:splitting1} is still nonlinear. Therefore, we consider a fractional step scheme
to solve \eqref{eqn:splitting1}
\begin{align*}
\frac{\mathbf{m}^* - \mathbf{m}^n}{\Delta t} &= \Delta_h\mathbf{m}^*\nonumber\\
\mathbf{m}^{n+1} &= \mathbf{m}^{n} - \mathbf{m}^{n}\times\mathbf{m}^*
\label{fractional}
\end{align*}
or
\begin{equation*}
\mathbf{m}^{n+1} = \mathbf{m}^{n} - \mathbf{m}^{n}\times(I-\Delta t
\Delta_h)^{-1}\mathbf{m}^n,
\end{equation*}
where $I$ is the identity matrix. This scheme is subject to strong stability constraint, and thus the implicit Gauss-Seidel scheme is introduced to overcome this issue.
Let
\begin{equation}
g_i^n = (I-\Delta t\Delta_h)^{-1}m_i^n,\quad i = 1,2,3.
\label{Thefirst_g}
\end{equation}
We then have
\begin{equation}
\begin{pmatrix}
m_1^{n+1}\\ m_2^{n+1}\\ m_3^{n+1}
\end{pmatrix}=
\begin{pmatrix}
m_1^n + (g_2^nm_3^n - g_3^nm_2^n)\\
m_2^n + (g_3^nm_1^{n+1} - g_1^{n+1}m_3^n)\\
m_3^n + (g_1^{n+1}m_2^{n+1} - g_2^{n+1}m_1^{n+1})
\end{pmatrix}.
\end{equation}
This scheme solve \eqref{eqn:splitting1} with unconditional stability.
\eqref{eqn:splitting2} is linear heat equation which can be solved easily.
However, the splitting scheme \eqref{eqn:splitting1} - \eqref{eqn:splitting2}
cannot preserve $\lvert\m\rvert = 1$, and thus a projection step needs to
be added.

For the full LLG equation \eqref{LLDimensionless}, the GSPM works as follows.
Define
\begin{equation}
\label{eqn:GSPM1}
\mathbf{h} = \epsilon\Delta\mathbf{m} + \mathbf{\hat f},
\end{equation}
where $\mathbf{\hat f}= -Q(m_2\mathbf{e}_2 + m_3\mathbf{e}_3) + \mathbf{h}_e + \mathbf{h}_s$.

The original GSPM \cite{NumGSPM2001} solves the equation \eqref{LLDimensionless} in three steps:
\begin{itemize}
\item Implicit Gauss-Seidel
\begin{align}
{g}_i^n &= (I - \Delta t\epsilon\Delta_h)^{-1}({m}_i^n + \Delta t{\hat f}_i^n),\ \ i = 2,3, \nonumber \\
{g}_i^* &= (I - \Delta t\epsilon\Delta_h)^{-1}({m}_i^* + \Delta t{\hat f}_i^n),\ \ i = 1,2,\label{eqn:GSPM2}
\end{align}
\begin{equation}
\label{eqn:GSPM3}
\begin{pmatrix}
{m}_1^* \\{m}_2^* \\{m}_3^*
\end{pmatrix} =
\begin{pmatrix}
{m}_1^n + ({g}_2^n{m}_3^n - {g}_3^n{m}_2^n) \\
{m}_2^n + ({g}_3^n{m}_1^* - {g}_1^*{m}_3^n) \\
{m}_3^n + ({g}_1^*{m}_2^* - {g}_2^*{m}_1^*)
\end{pmatrix}.
\end{equation}
\item Heat flow without constraints
\begin{equation}
\label{eqn:GSPM4}
\mathbf{\hat f}^* = -Q(m_2^*\mathbf{e}_2+m_3^*\mathbf{e}_3) + \mathbf{h}_e + \mathbf{h}_s^n,
\end{equation}
\begin{align}
\label{eqn:GSPM5}
\begin{pmatrix}
{m}_1^{**} \\{m}_2^{**} \\{m}_3^{**}
\end{pmatrix} =
\begin{pmatrix}
{m}_1^{*} + \alpha\Delta t(\epsilon\Delta_h{m}_1^{**} + \hat f_1^*) \\
{m}_2^{*} + \alpha\Delta t(\epsilon\Delta_h{m}_2^{**} + \hat f_2^*) \\
{m}_3^{*} + \alpha\Delta t(\epsilon\Delta_h{m}_3^{**} + \hat f_3^*)
\end{pmatrix}.
\end{align}
\item Projection onto $S^2$
\begin{align}
\label{eqn:GSPM6}
\begin{pmatrix}
{m}_1^{n+1} \\{m}_2^{n+1} \\{m}_3^{n+1}
\end{pmatrix} = \frac 1{|\mathbf{m}^{**}|}
\begin{pmatrix}
{m}_1^{**} \\{m}_2^{**} \\{m}_3^{**}
\end{pmatrix}.
\end{align}
\end{itemize}
Here the numerical stability of the original GSPM \cite{NumGSPM2001} was founded to be independent of gridsizes
but depend on the damping parameter $\alpha$. This issue was solved in \cite{ImproveGSPM2003} by replacing
\eqref{eqn:GSPM2} and \eqref{eqn:GSPM4} with
\[
{g}_i^* = (I - \Delta t\epsilon\Delta_h)^{-1}({m}_i^* + \Delta t{\hat f}_i^*),\ \ i = 1,2,
\]
and
\[
\mathbf{\hat f}^* = -Q(m_2^*\mathbf{e}_2+m_3^*\mathbf{e}_3) + \mathbf{h}_e + \mathbf{h}_s^*,
\]
respectively.
Update of the stray field is done using fast Fourier transform \cite{NumGSPM2001}.
It is easy to see that the GSPM solves $7$ linear systems of equations with constant coefficients
and updates the stray field using FFT $6$ times at each step.

\section{Two improved Gauss-Seidel projection methods for Landau-Lifshitz-Gilbert equation}\label{SectionImproved}

Based on the description of the original GSPM in \Cref{SectionGSPM}, we introduce two improved GSPMs
for LLG equation. The first improvement updates both the gyromagnetic term and the damping term
simultaneously, termed as Scheme A. The second improvement introduces two sets of approximate solution
with one set for implicit Gauss-Seidel step and the other set for projection in an alternating manner,
termed as Scheme B. Details are given in below.

\subsection{Scheme A}
The main improvement of Scheme A over the original GSPM is the combination of \eqref{eqn:GSPM1} - \eqref{eqn:GSPM5},
or \eqref{eqn:splitting1} - \eqref{eqn:splitting2}.
\begin{itemize}
\item Implicit-Gauss-Seidel
\begin{align}
\label{eqn:schemeA1}
{g}_i^n &= (I - \Delta t\epsilon\Delta_h)^{-1}({m}_i^n+\Delta t\hat{f}_i^n),\quad i = 1,2,3, \nonumber \\
{g}_i^{*} &= (I - \Delta t\epsilon\Delta_h)^{-1}({m}_i^{*}+\Delta t\hat{f}_i^*),\quad i = 1,2,
\end{align}
\begin{equation}
\label{eqn:schemeA2}
\begin{pmatrix}
{m}_1^* \\{m}_2^* \\{m}_3^*
\end{pmatrix} =
\begin{pmatrix}
m_1^n - (m_2^ng_3^n - m_3^ng_2^n) -
\alpha(m_1^ng_1^n+m_2^ng_2^n + m_3^ng_3^n)m_1^n + \alpha g_1^n \\
m_2^n - (m_3^ng_1^* - m_1^*g_3^n) -
\alpha(m_1^{*}g_1^{*}+m_2^ng_2^n + m_3^ng_3^n)m_2^n + \alpha g_2^n \\
m_3^n - (m_1^*g_2^* - m_2^*g_1^*) -
\alpha(m_1^{*}g_1^{*}+m_2^{*}g_2^{*} + m_3^ng_3^n)m_3^n + \alpha g_3^n
\end{pmatrix}.
\end{equation}
\item Projection onto $S^2$
\begin{align}
\label{eqn:schemeA3}
\begin{pmatrix}
{m}_1^{n+1} \\{m}_2^{n+1} \\{m}_3^{n+1}
\end{pmatrix} = \frac 1{|\mathbf{m}^{*}|}
\begin{pmatrix}
{m}_1^{*} \\{m}_2^{*} \\{m}_3^{*}
\end{pmatrix}.
\end{align}
\end{itemize}
It is easy to see that Scheme A solves $5$ linear systems of equations with constant coefficients
and uses FFT $5$ times at each step.

\subsection{Scheme B}

The main improvement of Scheme B over Scheme A is the introduction of two sets of approximate solutions,
one for \eqref{eqn:schemeA1} - \eqref{eqn:schemeA2} and the other for \eqref{eqn:schemeA3} and the update
of these two sets of solutions in an alternating manner.

Given the initialized $\mathbf{g}^0$
\begin{equation}
\label{eqn:schemeB1}
{g}_i^0 = (I - \Delta t\epsilon\Delta_h)^{-1}({m}_i^0+\Delta t\hat{f}_i^0),\ \ i = 1,2,3,
\end{equation}
Scheme B works as follows
\begin{itemize}
\item Implicit Gauss-Seidel
\begin{align}
\label{eqn:schemeB2}
{g}_i^{n+1} &= (I - \Delta t\epsilon\Delta_h)^{-1}({m}_i^{*}+\Delta t\hat{f}_i^{*}),\ \ i = 1,2,3
\end{align}
\begin{align}
\label{eqn:schemeB3}
m_1^{*} &= m_1^n - (m_2^ng_3^n - m_3^ng_2^n) -
\alpha(m_1^ng_1^n+m_2^ng_2^n + m_3^ng_3^n)m_1^n + \nonumber\\
&\ \ \ \alpha((m_1^n)^2 + (m_2^n)^2 + (m_3^n)^2)g_1^n\nonumber\\
m_2^{*} &= m_2^n - (m_3^ng_1^{n+1} - m_1^{*}g_3^n) -
\alpha(m_1^{*}g_1^{n+1}+m_2^ng_2^n + m_3^ng_3^n)m_2^n + \nonumber\\
&\ \ \ \alpha((m_1^{*})^2 + (m_2^n)^2 + (m_3^n)^2)g_2^n \nonumber\\
m_3^{*} &= m_3^n - (m_1^{*}g_2^{n+1} - m_2^{*}g_1^{n+1}) -
\alpha(m_1^{*}g_1^{n+1}+m_2^{*}g_2^{n+1} + m_3^ng_3^n)m_3^n + \nonumber\\
&\ \ \ \alpha((m_1^{*})^2 + (m_2^{*})^2 + (m_3^n)^2)g_3^n
\end{align}
\item Projection onto $S^2$
\begin{align}
\label{eqn:schemeB3}
\begin{pmatrix}
{m}_1^{n+1} \\{m}_2^{n+1} \\{m}_3^{n+1}
\end{pmatrix} = \frac 1{|\mathbf{m}^{*}|}
\begin{pmatrix}
{m}_1^{*} \\{m}_2^{*} \\{m}_3^{*}
\end{pmatrix}.
\end{align}
\end{itemize}
Here one set of approximate solution $\{\m^*\}$ is updated in the implicit Gauss-Seidel step
and the other set of approximate solution $\{\m^{n+1}\}$ is updated in the projection step.
Note that \eqref{eqn:schemeB2} is defined only for $\{\m^*\}$ which can be used in two successive temporal steps, and thus only $3$ linear systems of equations
with constant coefficients are solved at each step and $3$ FFT executions are used for
the stray field. The length of magnetization can be preserved in the time evolution.

The computational cost of GSPM and its improvements comes from solving
the linear systems of equations with constant coefficients.
To summarize, we list the number of linear systems of equations to be solved and
the number of FFT executions to be used at each step for the original GSPM \cite{ImproveGSPM2003},
Scheme A, and Scheme B in  \cref{savetime}.
The savings represent the ratio between costs of two improved schemes over that of the original GSPM.

\begin{table}[ht]
	\centering
	\begin{tabular}{c|c|c|c|c}
		\hline
		GSPM Scheme & Number of linear systems & Saving & Execution of FFT & Saving \\
		\hline
		Original & $7$ & $0$   & $4$ & $0$   \\
		Scheme A & $5$ & $2/7$ & $3$ & $1/4$ \\
		Scheme B & $3$ & $4/7$ & $3$ & $1/4$ \\
		\hline
	\end{tabular}
	\caption{The number of linear systems of equations to be solved and the number of FFT executions to be used
		at each step for the original GSPM \cite{ImproveGSPM2003}, Scheme A, and Scheme B.
		The savings represent the ratio between costs of two improved schemes over that of the original GSPM.}
	\label{savetime}
\end{table}

\section{Numerical Experiments}\label{sectionNumerical}

In this section, we compare the original GSPM \cite{NumGSPM2001, ImproveGSPM2003}, Scheme A, and Scheme B
via a series of examples in both 1D and 3D, including accuracy check and efficiency check, unconditional stability
with respect to both the grid size and the damping parameter, hysteresis loops, and magnetization profiles.
For convenience, we define
\[
\mathrm{ratio}-i = \frac{\mathrm{Time(GSPM)} - \mathrm{Time(Scheme}\; i)}{\mathrm{Time(GSPM)}},
\]
for $i =$ A and B, which quantifies the improved efficiency of Scheme A and Scheme B over the original GSPM \cite{NumGSPM2001, ImproveGSPM2003}.

\subsection{Accuracy Test}

\begin{example}[1D case]

In 1D, we choose the exact solution over the unit interval $\Omega=[0,1]$
\[
\mathbf{m}_e = (\cos(\bar x)\sin(t), \sin(\bar x)\sin(t), \cos(t)),
\]
which satisfies
\[
\m_t = - \m\times \m_{xx} - \alpha\m\times(\m\times\m_{xx}) + \mathbf{f}
\]
with $\bar x=x^2(1-x)^2$, and $\mathbf{f}=\mathbf{m}_{et} +
\mathbf{m}_e\times\mathbf{m}_{exx} + \alpha\mathbf{m}_e\times(\mathbf{m}_e\times\mathbf{m}_{exx})$.
Parameters are $\alpha=0.00001$ and $T=5.0e-2$.

We first show the error $\|\mathbf{m}_e - \mathbf{m}_h\|_{\infty}$ with $\mathbf{m}_h$ being the numerical solution
with respect to the temporal stepsize $\Delta t$ and the spatial stepsize $\Delta x$.
As shown in \cref{1Dfigure}(a) and \cref{1Dfigure}(c), suggested by the least squares fitting,
both first-order accuracy in time and second-order accuracy in space
are observed. Meanwhile, we record the CPU time as a function of accuracy (error) by varying
the temporal stepsize and the spatial stepsize in \cref{1Dfigure}(b) and \cref{1Dfigure}(d),
\cref{1DTemporaltable} and \cref{1DSpatialtable}, respectively. In addition, from \cref{1DTemporaltable}
and \cref{1DSpatialtable}, the saving of Scheme A over GSPM is about $2/7$, which equals $1-5/7$,
and the saving of Scheme B over GSPM is about $4/7$, respectively. This observation is in good agreement
with the number of linear systems being solved at each step for these three methods, as shown in \cref{savetime}.
\begin{figure}[ht]
\centering
\subfloat[Temporal accuracy]
{\label{temporal1D}\includegraphics[width=2.5in]{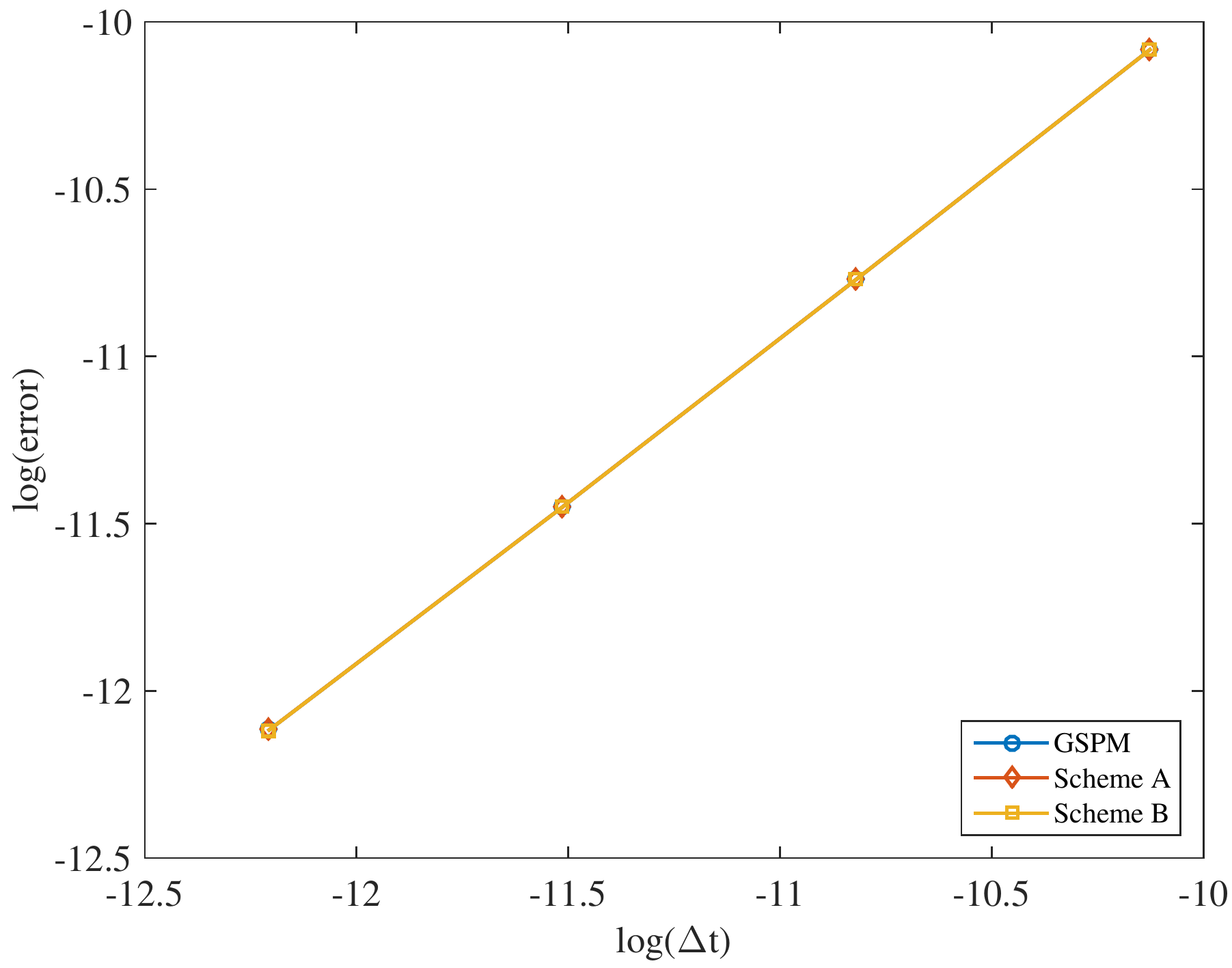}}
\subfloat[CPU time versus approximation error ($\Delta t$)]
{\label{temporalTime1D}\includegraphics[width=2.5in]{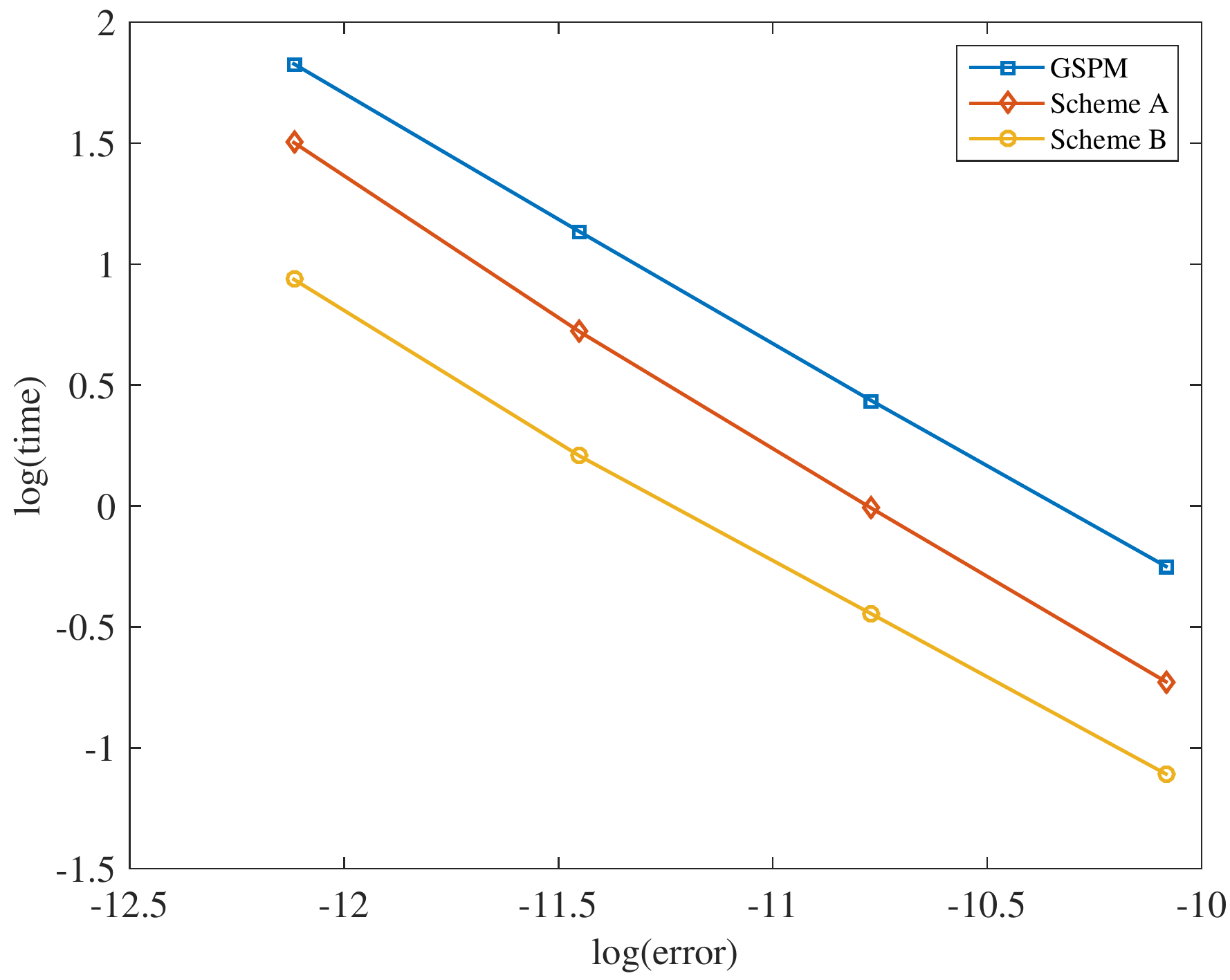}}
\quad
\subfloat[Spatial accuracy]
{\label{spatial1D}\includegraphics[width=2.5in]{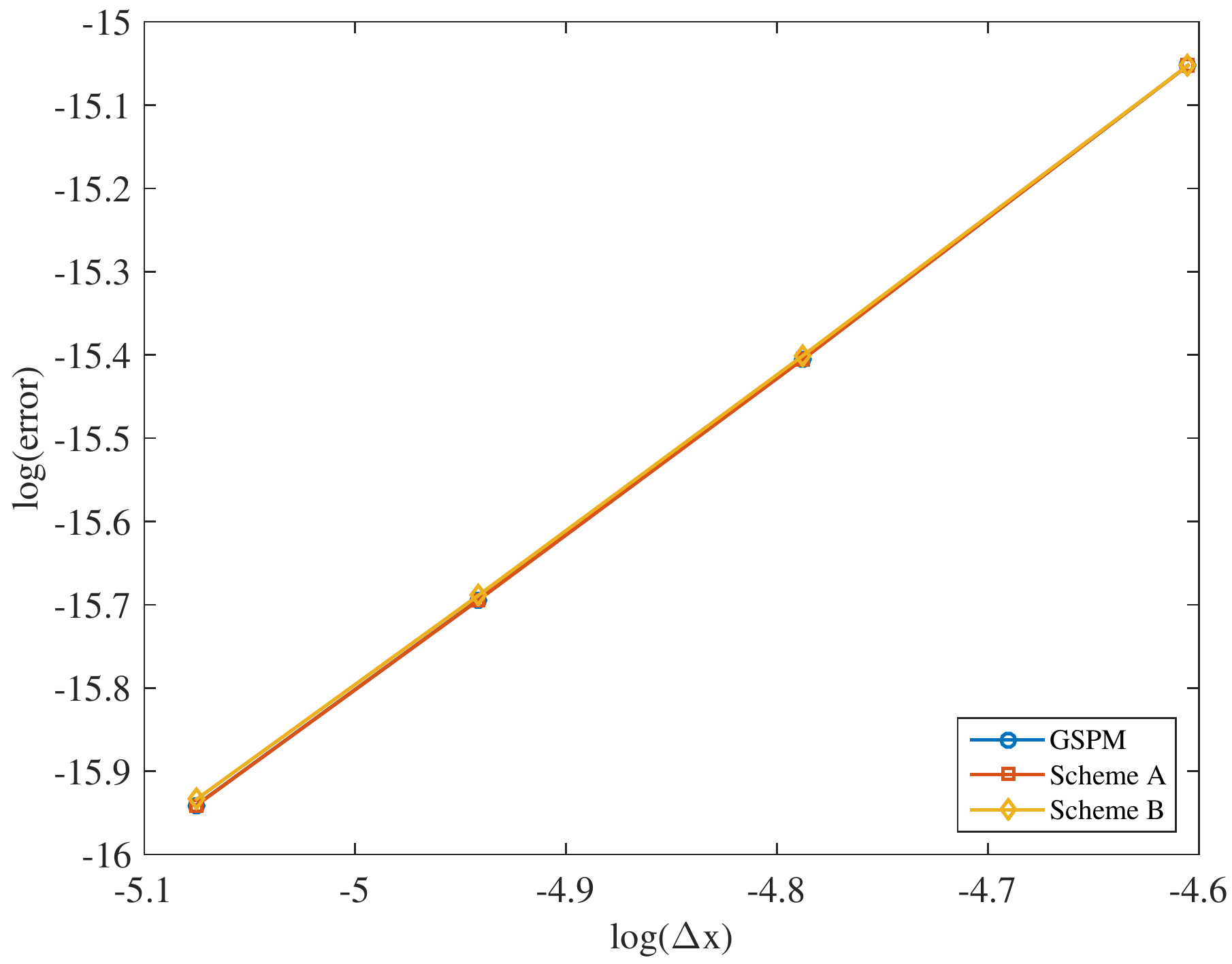}}
\subfloat[CPU time versus approximation error ($\Delta x$)]
{\label{spatialTime1D}\includegraphics[width=2.5in]{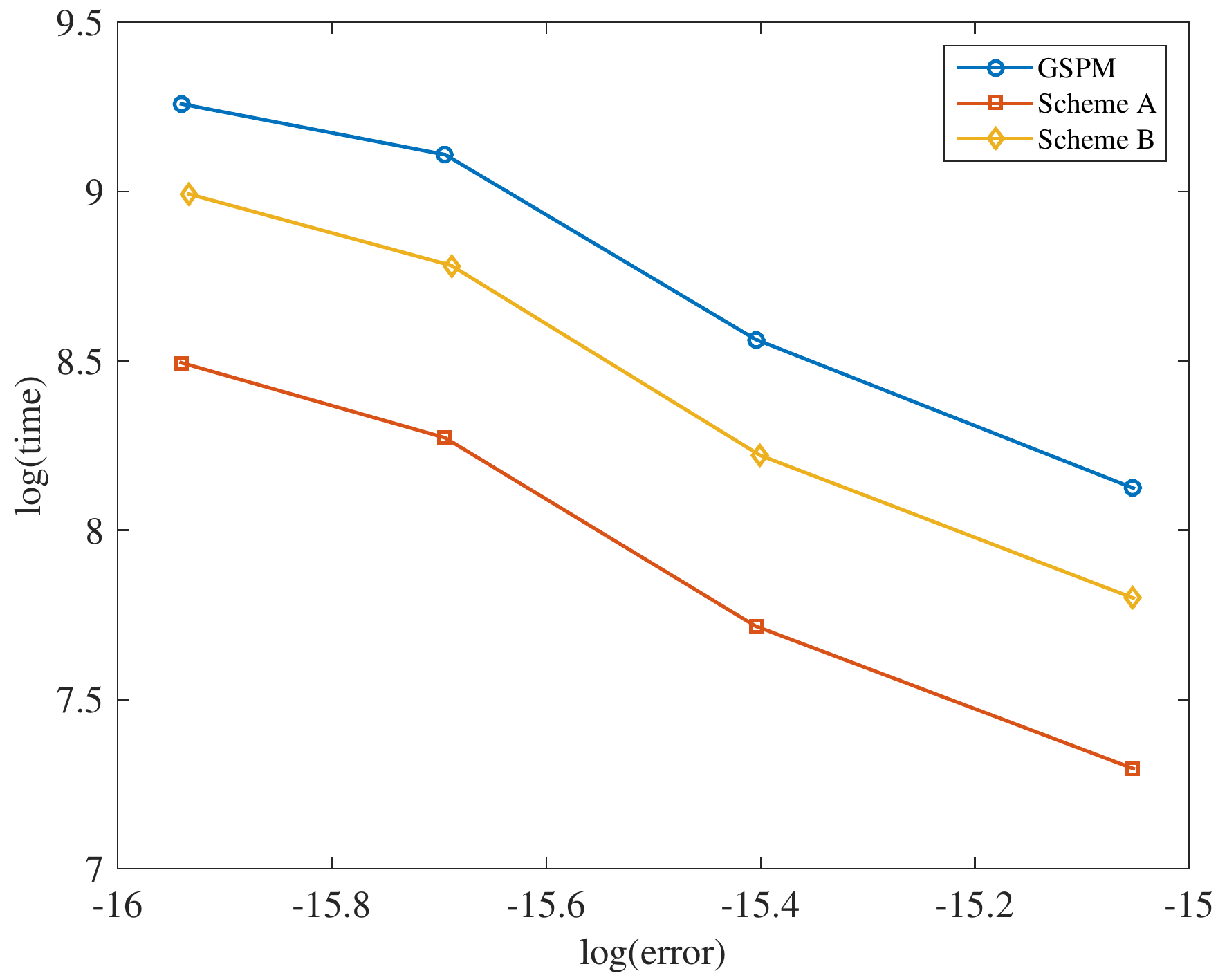}}
\caption{Approximation error and CPU time in 1D. (a) Approximation error as a function of the temporal step size;
(b) CPU time as a function of the approximation error when $\Delta t$ is varied and $\Delta x$ is fixed;
(c) Approximation error as a function of the spatial step size;
(d) CPU time as a function of the approximation error when $\Delta x$ is varied and $\Delta t$ is fixed.}
\label{1Dfigure}
\end{figure}

\begin{table}[ht]
\centering
\begin{tabular}{l|c|c|c|c|c}
\hline
\backslashbox{CPU time}{$\Delta t$} & T/1250 & T/2500 & T/5000 & T/10000 & Reference\\
\hline
GSPM & 7.7882e-01 & 1.5445e+00 & 3.1041e+00 & 6.2196e+00 & -\\
Scheme A & 4.8340e-01 & 9.9000e-01 & 2.0527e+00 & 4.4917e+00 & -\\
Scheme B & 3.3010e-01 & 6.3969e-01 & 1.2281e+00 & 2.5510e+00 & -\\
\hline
ratio-A & 0.38 & 0.36 & 0.34 & 0.28 & 0.29(2/7)\\
ratio-B & 0.58 & 0.59 & 0.60 & 0.59 & 0.57(4/7)\\
\hline
\end{tabular}
\caption{Recorded CPU time in 1D with respect to the approximation error when only $\Delta t$ is varied and $\Delta x = 1/100$.}
\label{1DTemporaltable}
\end{table}

\begin{table}[ht]
\centering
\begin{tabular}{l|c|c|c|c|c}
\hline
\backslashbox{CPU time}{$\Delta x$} & 1/100 & 1/120 & 1/140 & 1/160 & Reference\\
\hline
GSPM & 3.3752e+03 & 5.2340e+03 & 9.0334e+03 & 1.0495e+04 & -\\
Scheme A & 2.4391e+03 & 3.7175e+03 & 6.5149e+03 & 8.0429e+03 & -\\
Scheme B & 1.4740e+03 & 2.2448e+03 & 3.9152e+03 & 4.8873e+03 & -\\
\hline
ratio-A & 0.28 & 0.29 & 0.28 & 0.23 & 0.29(2/7)\\
ratio-B & 0.56 & 0.57 & 0.57 & 0.53 & 0.57(4/7)\\
\hline
\end{tabular}
\caption{Recorded CPU time in 1D with respect to the approximation error when only $\Delta x$ is varied and $\Delta t = 1.0e-8$.}
\label{1DSpatialtable}
\end{table}

\end{example}

\begin{example}[3D case]

In 3D, we choose the exact solution over $\Omega=[0,2]\times[0,1]\times[0,0.2]$
\[
\mathbf{m}_e = (\cos(\bar x \bar y \bar z)\sin(t),\sin(\bar x \bar y \bar z)\sin(t), \cos(t)),
\]
which satisfies
\[
\m_t=-\m\times\Delta\m - \alpha\m\times(\m\times\Delta \m) + \mathbf{f}
\]
with $\bar x=x^2(1-x)^2$, $\bar y=y^2(1-y)^2$, $\bar z=z^2(1-z)^2$ and $\mathbf{f}=\mathbf{m}_{et} +
\mathbf{m}_e\times\Delta\mathbf{m}_e + \alpha\mathbf{m}_e\times(\mathbf{m}_e\times\Delta\mathbf{m}_e)$.
Parameters are $T = 1.0e-05$ and $\alpha = 0.01$.

Like in the 1D case, we first show the error $\|\mathbf{m}_e - \mathbf{m}_h\|_{\infty}$ with $\mathbf{m}_h$ being the numerical solution
with respect to the temporal stepsize $\Delta t$ and the spatial stepsize $\Delta x$.
As shown in \cref{3Dfigure}(a) and \cref{3Dfigure}(c), suggested by the least squares fitting,
both first-order accuracy in time and second-order accuracy in space
are observed. Meanwhile, we record the CPU time as a function of accuracy (error) by varying
the temporal stepsize and the spatial stepsize in \cref{3Dfigure}(b) and \cref{3Dfigure}(d),
\cref{3DTemporaltable} and \cref{3DSpatialtable}, respectively. In addition, from \cref{3DTemporaltable}
and \cref{3DSpatialtable}, the saving of Scheme A over GSPM is about $2/7$,
and the saving of Scheme B over GSPM is about $4/7$, respectively. This observation is in good agreement
with the number of linear systems being solved at each step for these three methods, as shown in \cref{savetime}.

\begin{figure}[ht]
	\centering
	\subfloat[Temporal accuracy]
	{\label{temporal3D}\includegraphics[width=6cm]{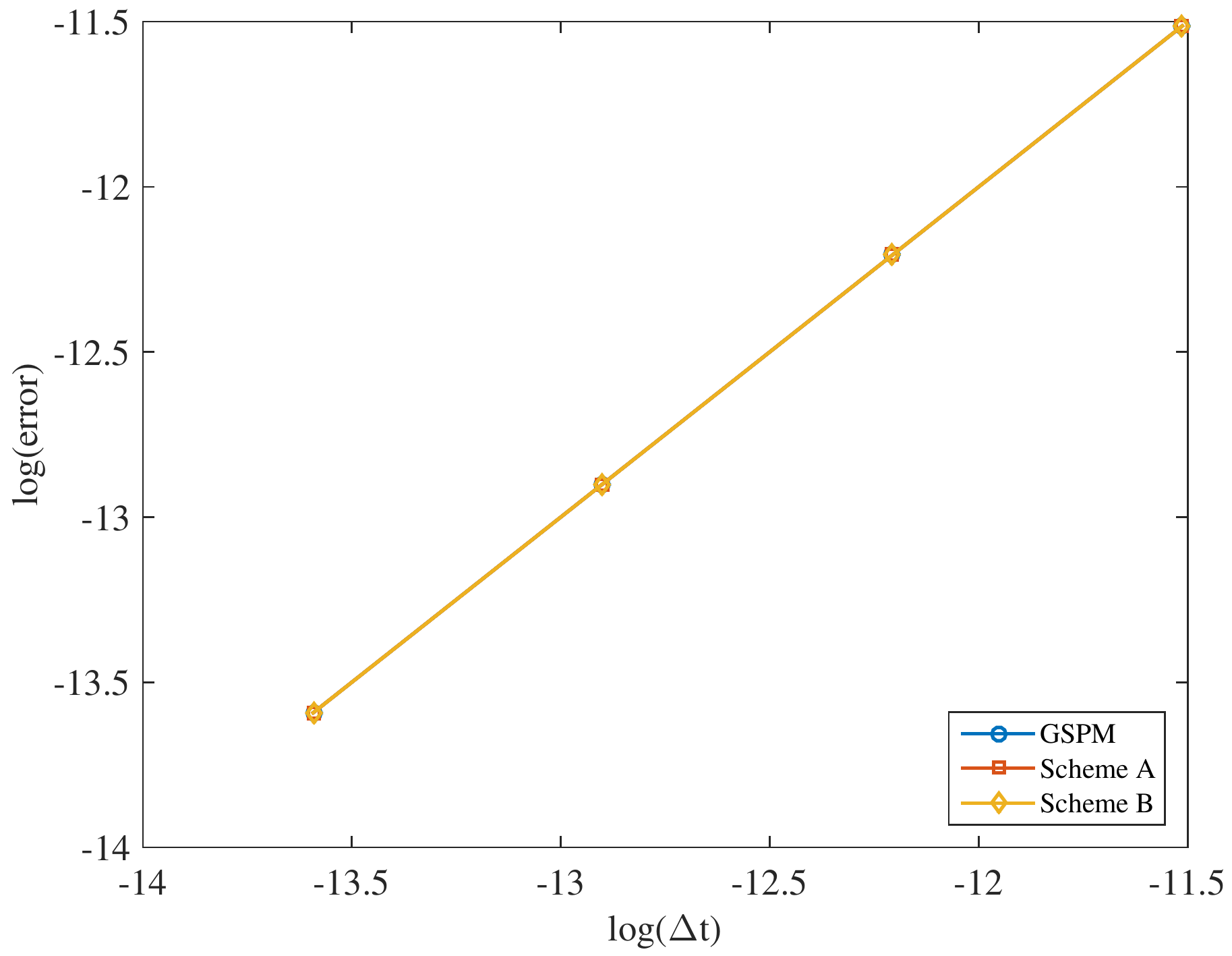}}
	\subfloat[CPU time versus approximation error ($\Delta t$)]
	{\label{temporalTime3D}\includegraphics[width=6cm]{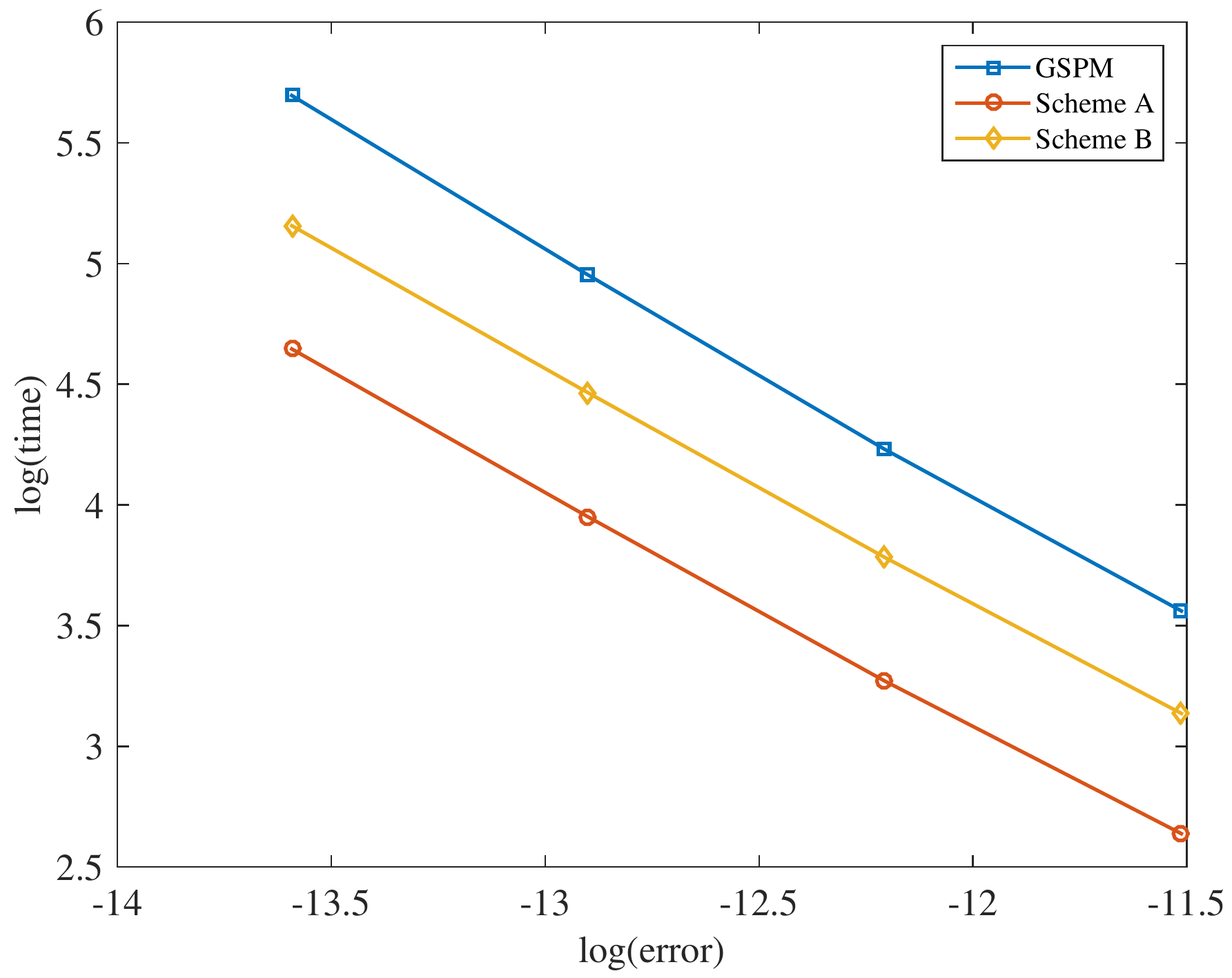}}
	\quad
	\subfloat[Spatial accuracy]
	{\label{spatial3D}\includegraphics[width=6cm]{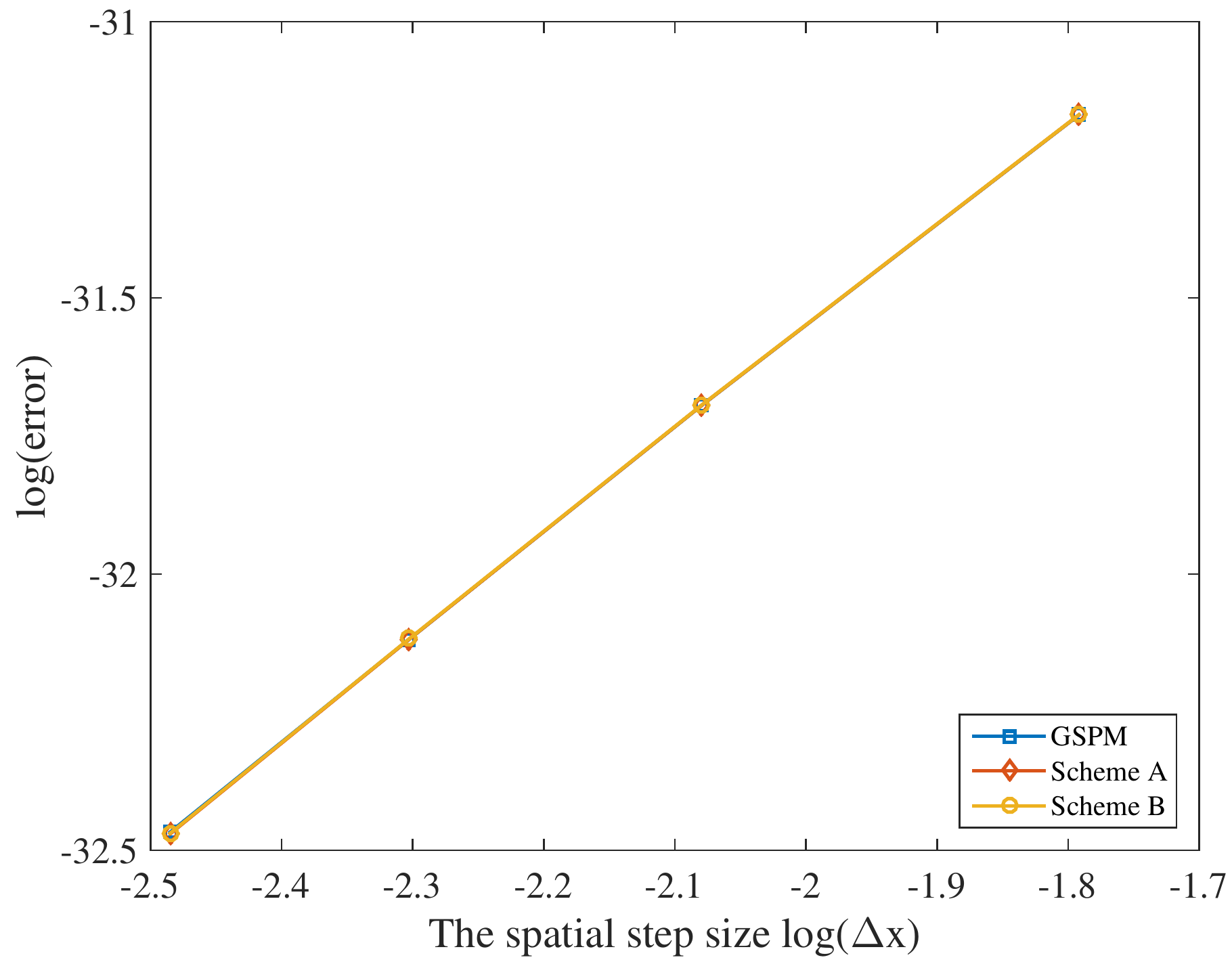}}
	\subfloat[CPU time versus approximation error ($\Delta x$)]
	{\label{spatialTime3D}\includegraphics[width=6cm]{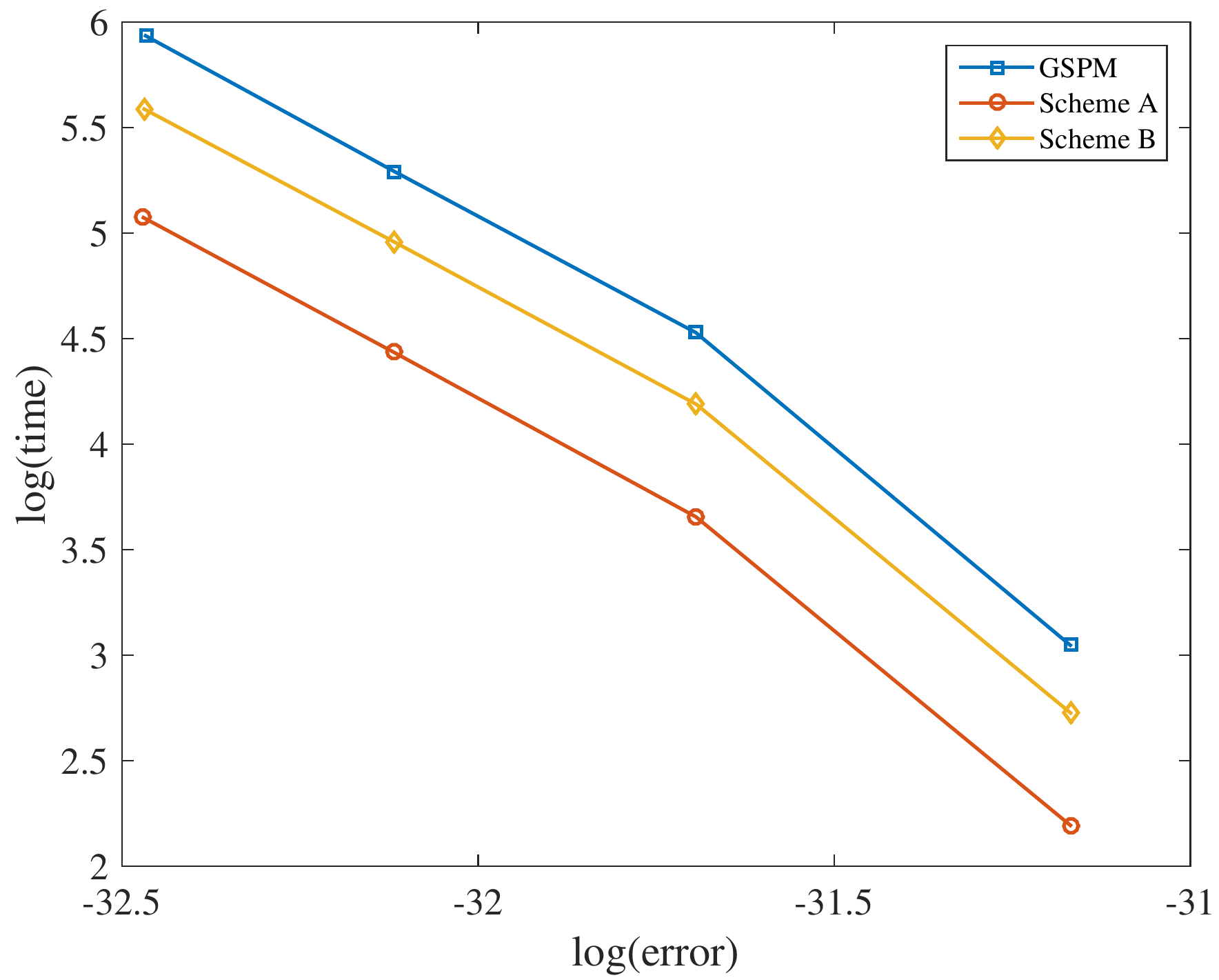}}
	\caption{Approximation error and CPU time in 3D. (a) Approximation error as a function of the temporal step size;
		(b) CPU time as a function of the approximation error when $\Delta t$ is varied and $\Delta x = \Delta y = \Delta z$ is fixed;
		(c) Approximation error as a function of the spatial step size;
		(d) CPU time as a function of the approximation error when space is varied uniformly and $\Delta t$ is fixed.}
	\label{3Dfigure}
\end{figure}

\begin{table}[ht]
\centering
\begin{tabular}{l|c|c|c|c|c}
\hline
\backslashbox{CPU time}{$\Delta t$} & T/10 & T/20 & T/40 & T/80 & Reference\\
\hline
GSPM & 3.5188e+01 & 6.8711e+01 & 1.4146e+02 & 2.9769e+02 & -\\
Scheme A & 2.3015e+01 & 4.3920e+01 & 8.6831e+01 & 1.7359e+02 &-\\
Scheme B & 1.3984e+01 & 2.6313e+01 & 5.1928e+01 & 1.0415e+02 & -\\
\hline
ratio-A & 0.35 & 0.36 & 0.39 & 0.42 & 0.29(2/7)\\
ratio-B & 0.60 & 0.62 & 0.63 & 0.65 & 0.57(4/7)\\
\hline
\end{tabular}
\caption{Recorded CPU time in 3D with respect to the approximation error when only $\Delta t$ is varied and the spatial mesh is $128\times64\times10$.}
\label{3DTemporaltable}
\end{table}

\begin{table}[ht]
\centering
\begin{tabular}{l|c|c|c|c|c}
\hline
\backslashbox{CPU time}{$\Delta x$}  & 1/6 & 1/8 & 1/10 & 1/12 & Reference\\
\hline
GSPM & 2.1066e+01 & 9.2615e+01 & 1.9879e+02 & 3.7820e+02 & -\\
Scheme A & 1.5278e+01 & 6.5953e+01 & 1.4215e+02 & 2.6725e+02 & -\\
Scheme B & 8.9698e+00 & 3.8684e+01 & 8.4291e+01 & 1.5977e+02 & -\\
\hline
ratio-A & 0.27 & 0.29 & 0.28 & 0.29 & 0.29(2/7)\\
ratio-B & 0.57 & 0.58 & 0.58 & 0.58 & 0.57(4/7)\\
\hline
\end{tabular}
\caption{Recorded CPU time in 3D with respect to the approximation error when only the spatial gridsize is varied with $\Delta x = \Delta y = \Delta z$ and $\Delta t = 1.0e-09$.}
\label{3DSpatialtable}
\end{table}

\end{example}

It worths mentioning that all these three methods are tested to be unconditionally stable with respect to the spatial gridsize and the temporal stepsize.

\subsection{Micromagnetic Simulations}

To compare the performance of Scheme A and Scheme B with GSPM, we have
carried out micromagnetic simulations of the full LLG equation with realistic material parameters.
In all our following simulations, we consider a thin film ferromagnet of size $\Omega=1\;\mu\textrm{m}\times1\;\mu\textrm{m}\times0.02\;\mu\textrm{m}$
with the spatial gridsize $4\;\textrm{nm}\times4\;\textrm{nm}\times4\;\textrm{nm}$ and the temporal stepsize $\dt=1$ picosecond.
The demagnetization field (stray field) is calculated via FFT \cite{NumGSPM2001,ImproveGSPM2003}.

\subsubsection{Comparison of hysteresis loops}
The hysteresis loop is calculated in the following way. First, a positive external field
$H_0=\mu_0 H$ is applied and the system is allowed to reach a stable state. Afterwards,
the external field is reduced by a certain amount and the system is relaxed to a stable state again.
The process continues until the external field attains a negative field of strength $H_0$.
Then the external field starts to increase and the system relaxes until the initial applied external
field $H_0$ is approached. In the hysteresis loop, we can monitor the magnetization dynamics and
plot the average magnetization at the stable state as a function of the strength of the external field.
The stopping criterion for a steady state is that the relative change of the total energy is less than $10^{-7}$.
The applied field is parallel to the $x$ axis. The initial state we take is the uniform state and the damping parameter $\alpha=0.1$.

\begin{figure}[ht]
\centering
\includegraphics[width=4in]{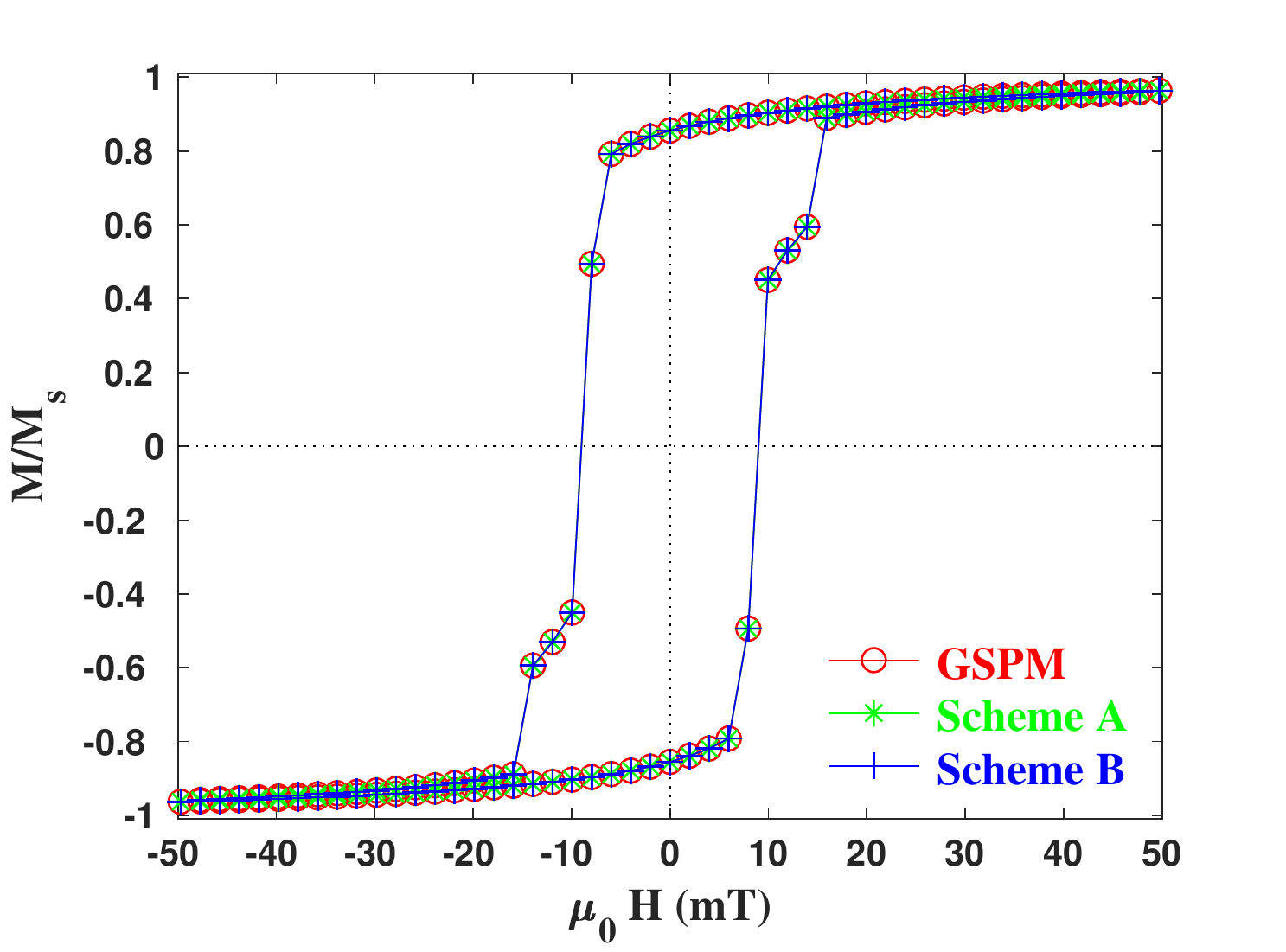}
\caption{\footnotesize Comparison of hysteresis loops for GSPM, Scheme A and Scheme B. Profiles of the average magnetization
of these three methods are in quantitative agreements with approximately the same switch field $9\;(\pm 0.4)\;\mathrm{mT}$.
The applied field is parallel to the $x$ axis and the initial state is the uniform state.}\label{loop}
\end{figure}

In \cref{loop}, we compare the average magnetization in the hysteresis loop simulated by GSPM, Scheme A and
Scheme B. Profiles of the average magnetization of these three methods are in quantitative agreements with
approximately the same switch field $9\;(\pm 0.4)\;\mathrm{mT}$.

\subsubsection{Comparison of magnetization profiles}

\begin{figure}[htbp]
	\centering
	\subfloat[Angle profile ($\alpha=0.1$)]{\label{GSPM_angle_1}\includegraphics[width=2.5in]{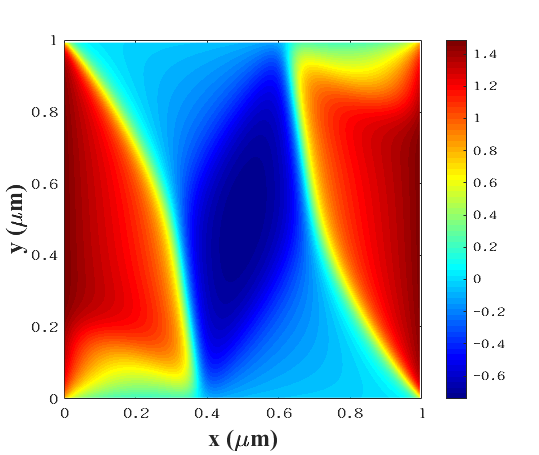}}
	\subfloat[Magnetization profile ($\alpha=0.1$)]{\label{GSPM_mag_1}\includegraphics[width=2.5in]{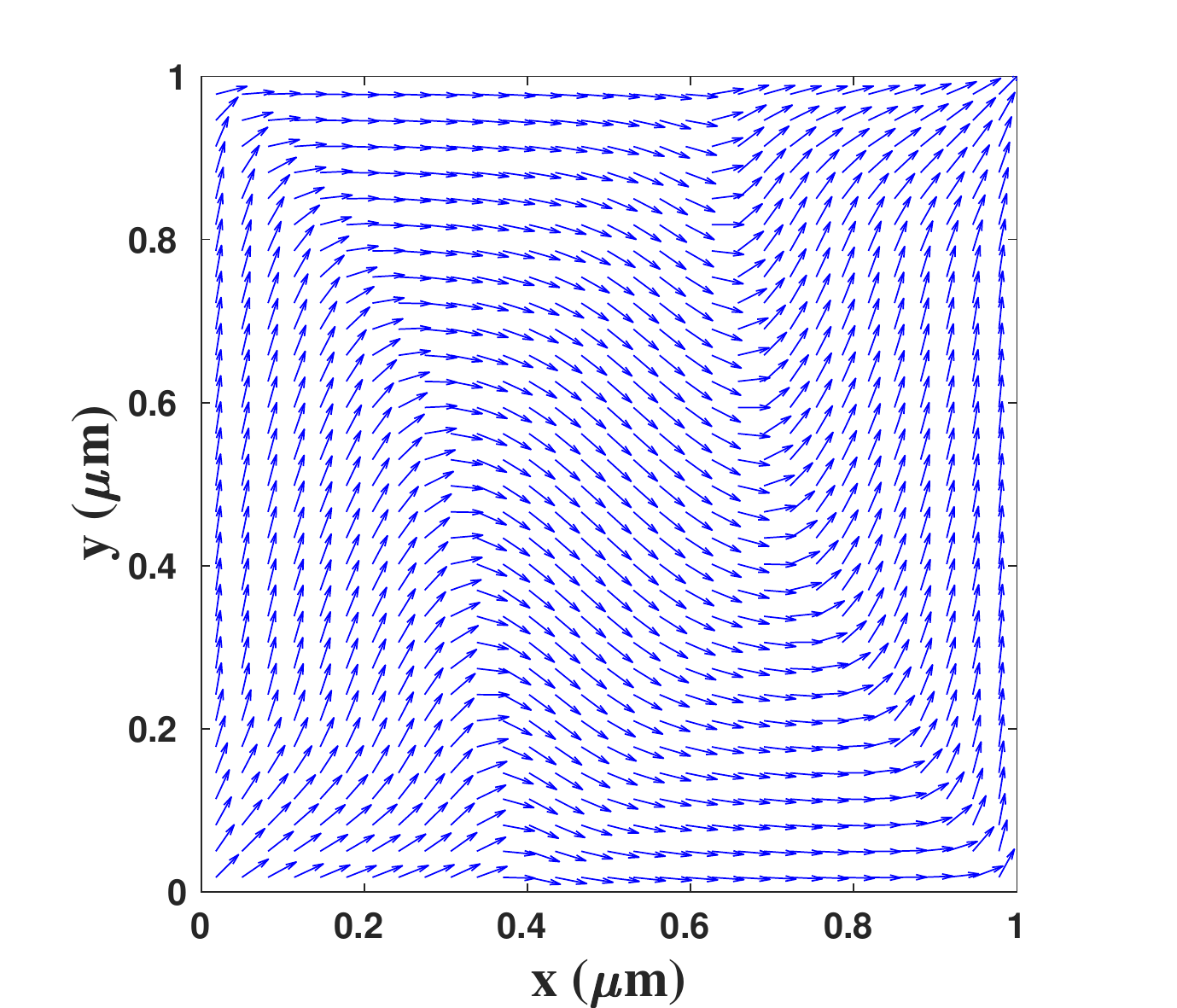}}
	\quad
	\subfloat[Angle profile ($\alpha=0.01$)]{\label{GSPM_angle_2}\includegraphics[width=2.5in]{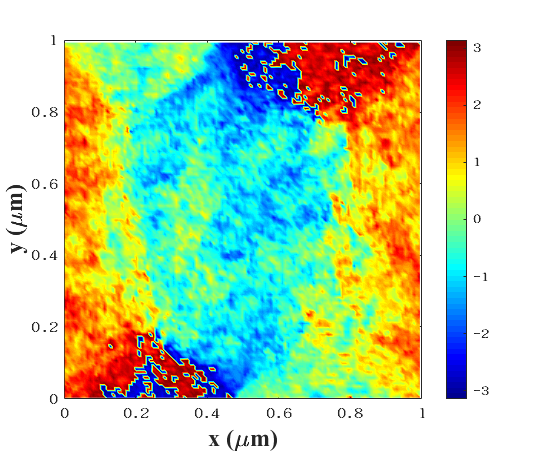}}
	\subfloat[Magnetization profile ($\alpha=0.01$)]{\label{GSPM_mag_2}\includegraphics[width=2.5in]{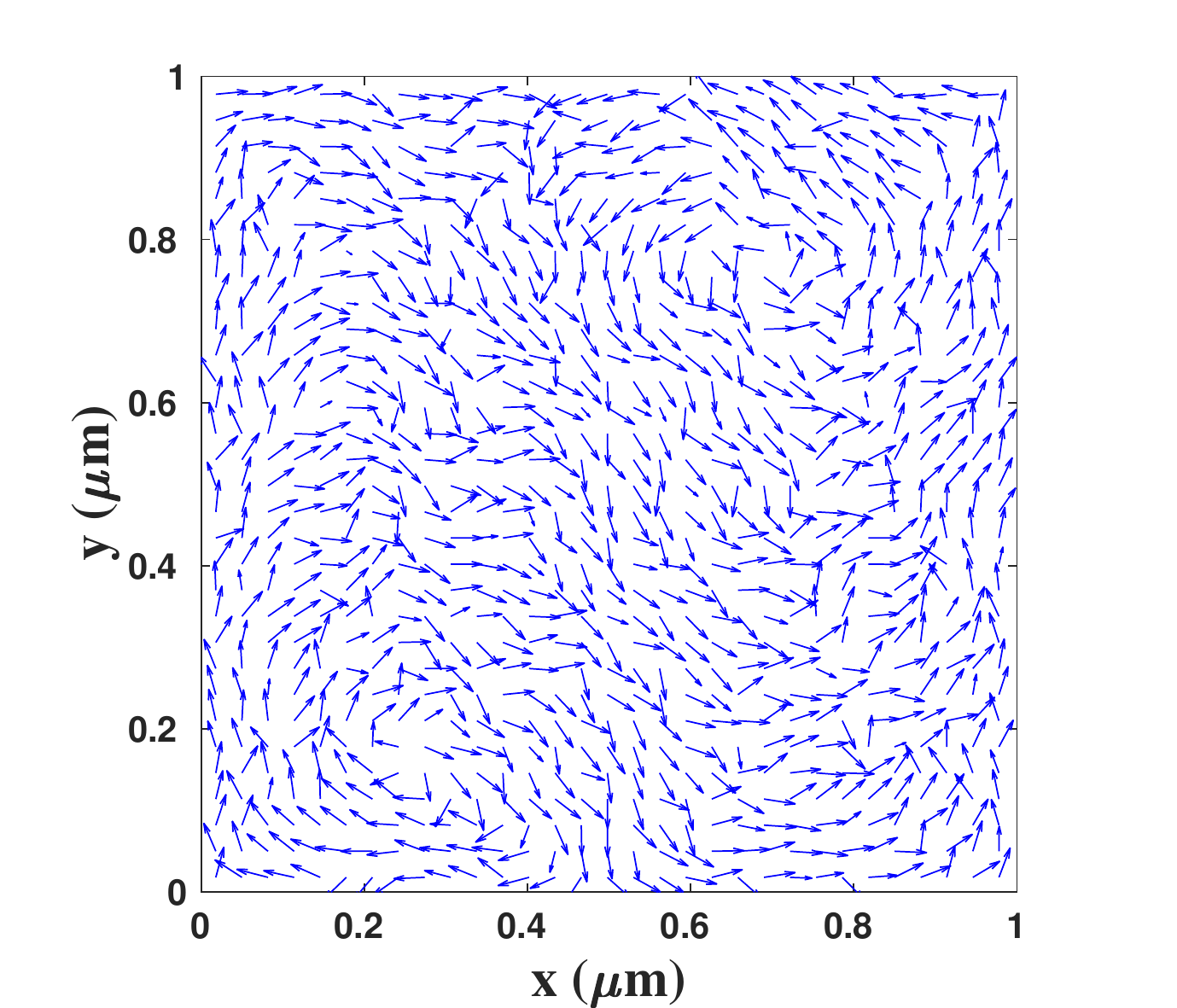}}
	\caption{Simulation of the full Landau-Lifshitz-Gilbert equation using GSPM without any external field. The magnetization on the centered slice of the material in the $xy$ plane is used.
Top row: $\alpha=0.1$; Bottom row: $\alpha=0.01$. Left column: a color plot of the angle between the in-plane magnetization and the $x$ axis;
Right column: an arrow plot of the in-plane magnetization.}\label{GSPM}
\end{figure}

\begin{figure}[htbp]
	\centering
	\subfloat[Angle profile ($\alpha=0.1$)]{\label{IGS_angle_1}\includegraphics[width=2.5in]{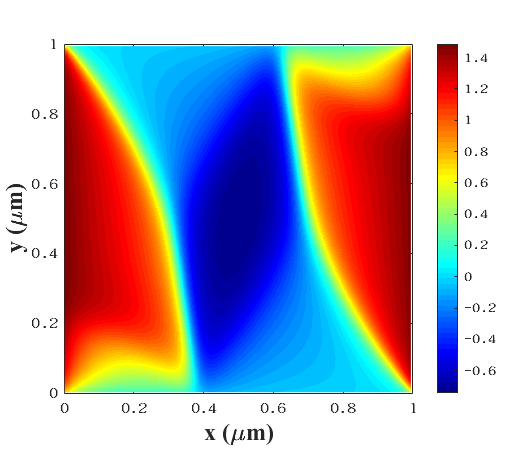}}
	\subfloat[Magnetization profile ($\alpha=0.1$)]{\label{IGS_mag_1}\includegraphics[width=2.5in]{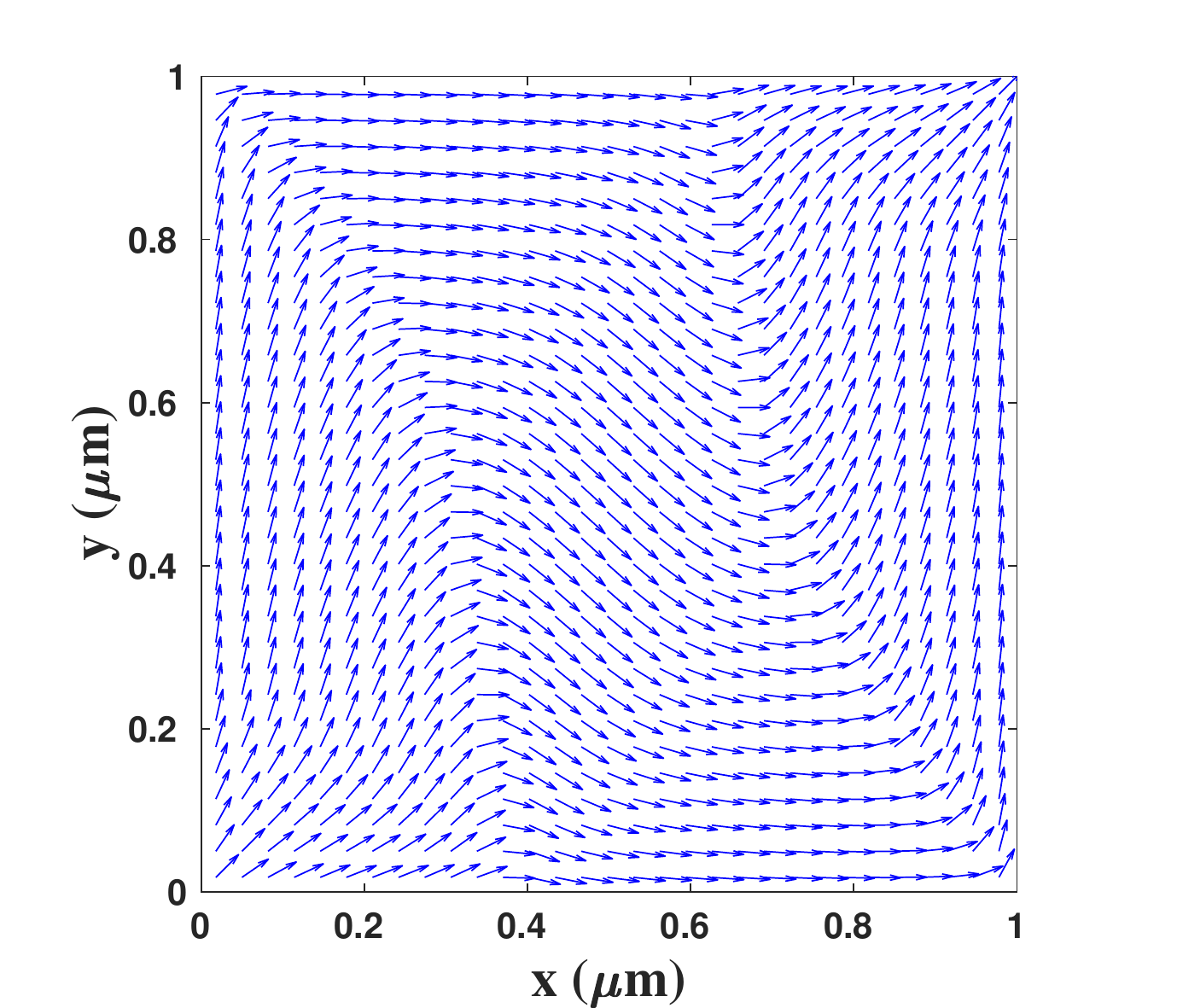}}
	\quad
	\subfloat[Angle profile ($\alpha=0.01$)]{\label{IGS_angle_2}\includegraphics[width=2.5in]{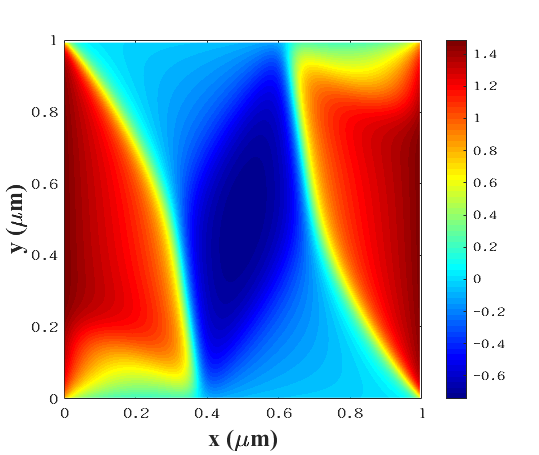}}
	\subfloat[Magnetization profile ($\alpha=0.01$)]{\label{IGS_mag_2}\includegraphics[width=2.5in]{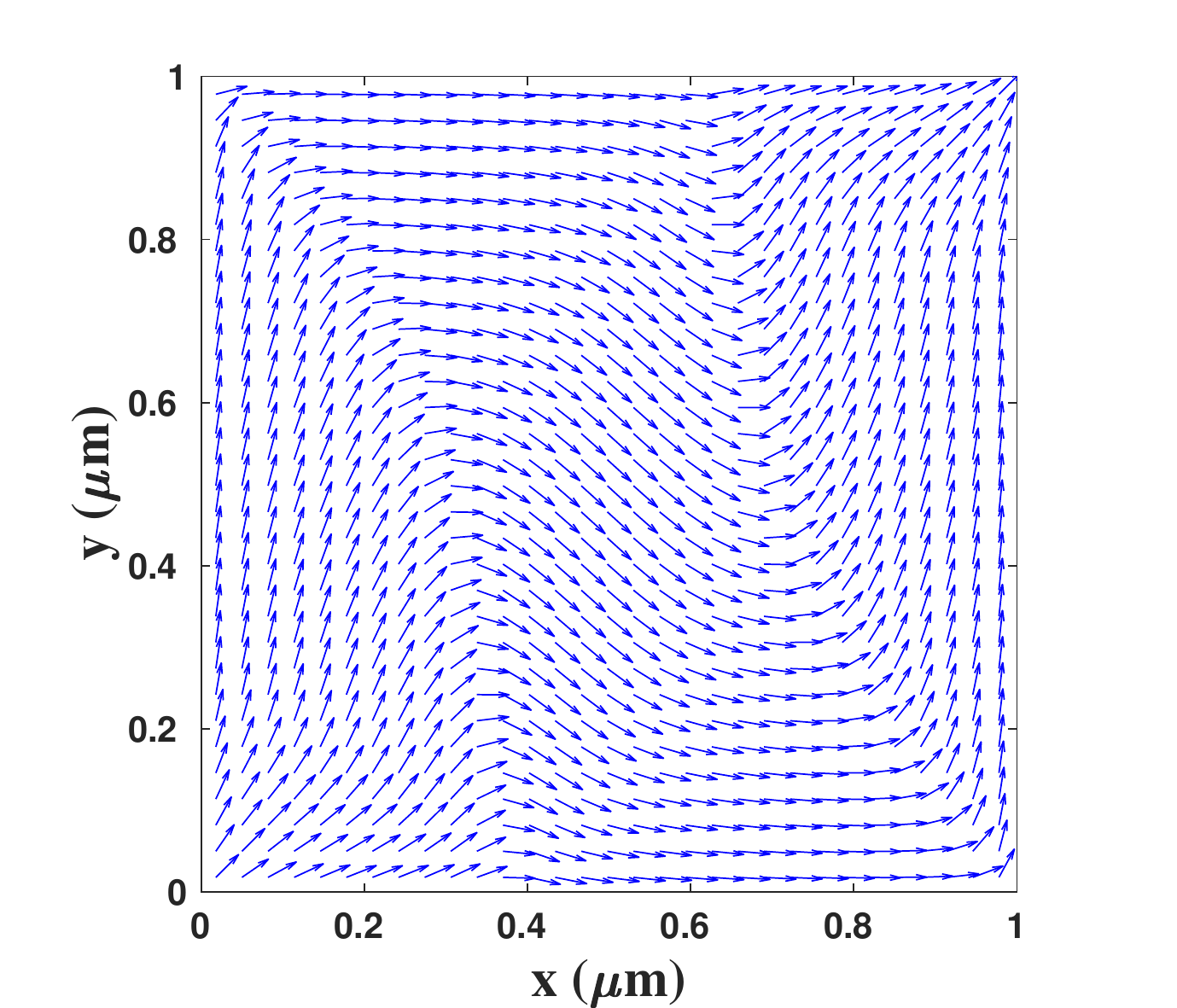}}
	\caption{Simulation of the full Landau-Lifshitz-Gilbert equation using Scheme A without any external field. The magnetization on the centered slice of the material in the $xy$ plane is used.
		Top row: $\alpha=0.1$; Bottom row: $\alpha=0.01$. Left column: a color plot of the angle between the in-plane magnetization and the $x$ axis;
		Right column: an arrow plot of the in-plane magnetization.}
	\label{IGS}
\end{figure}

\begin{figure}[htbp]
	\centering
	\subfloat[Angle profile ($\alpha=0.1$)]{\label{IGSP_angle_1}\includegraphics[width=2.5in]{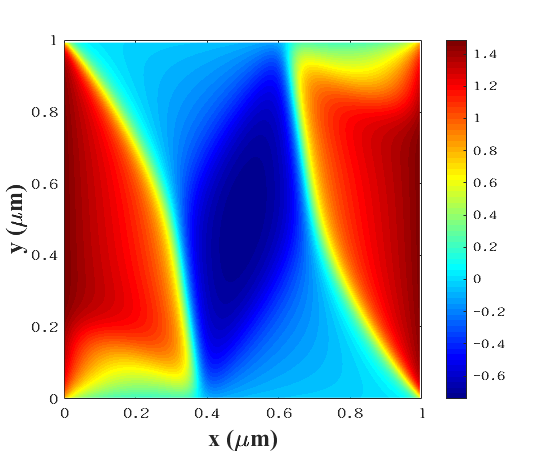}}
	\subfloat[Magnetization profile ($\alpha=0.1$)]{\label{IGSP_mag_1}\includegraphics[width=2.5in]{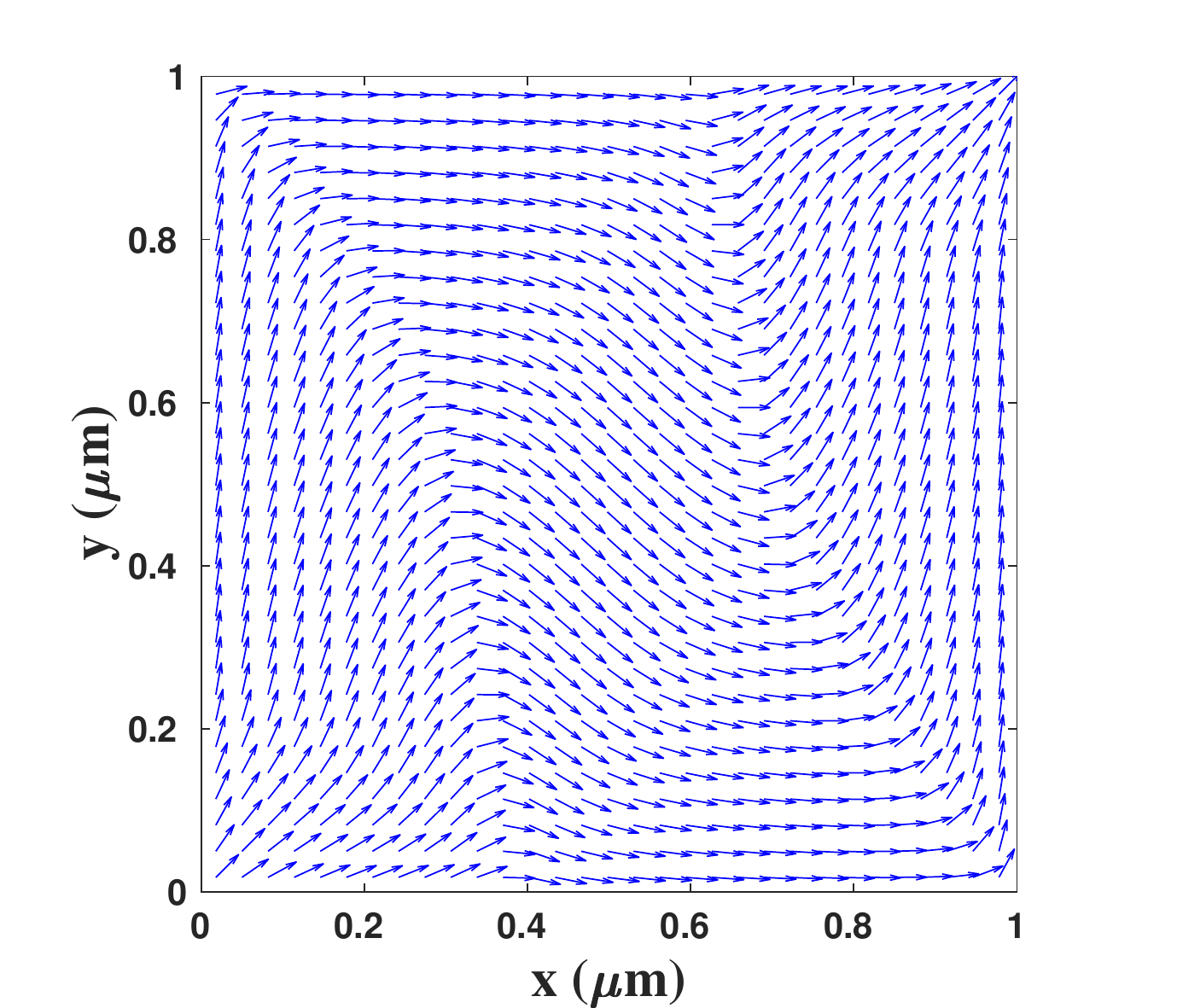}}
	\quad
	\subfloat[Angle profile ($\alpha=0.01$)]{\label{IGSP_angle_2}\includegraphics[width=2.5in]{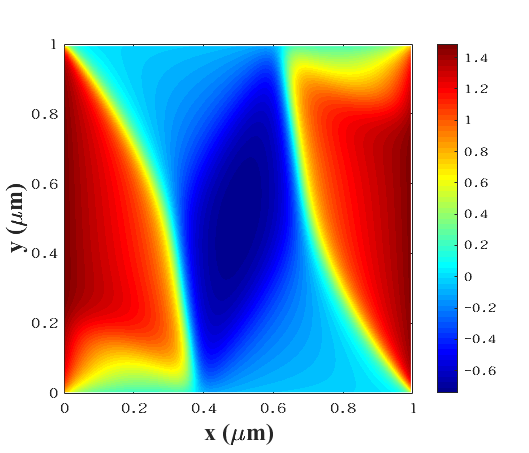}}
	\subfloat[Magnetization profile ($\alpha=0.01$)]{\label{IGSP_mag_2}\includegraphics[width=2.5in]{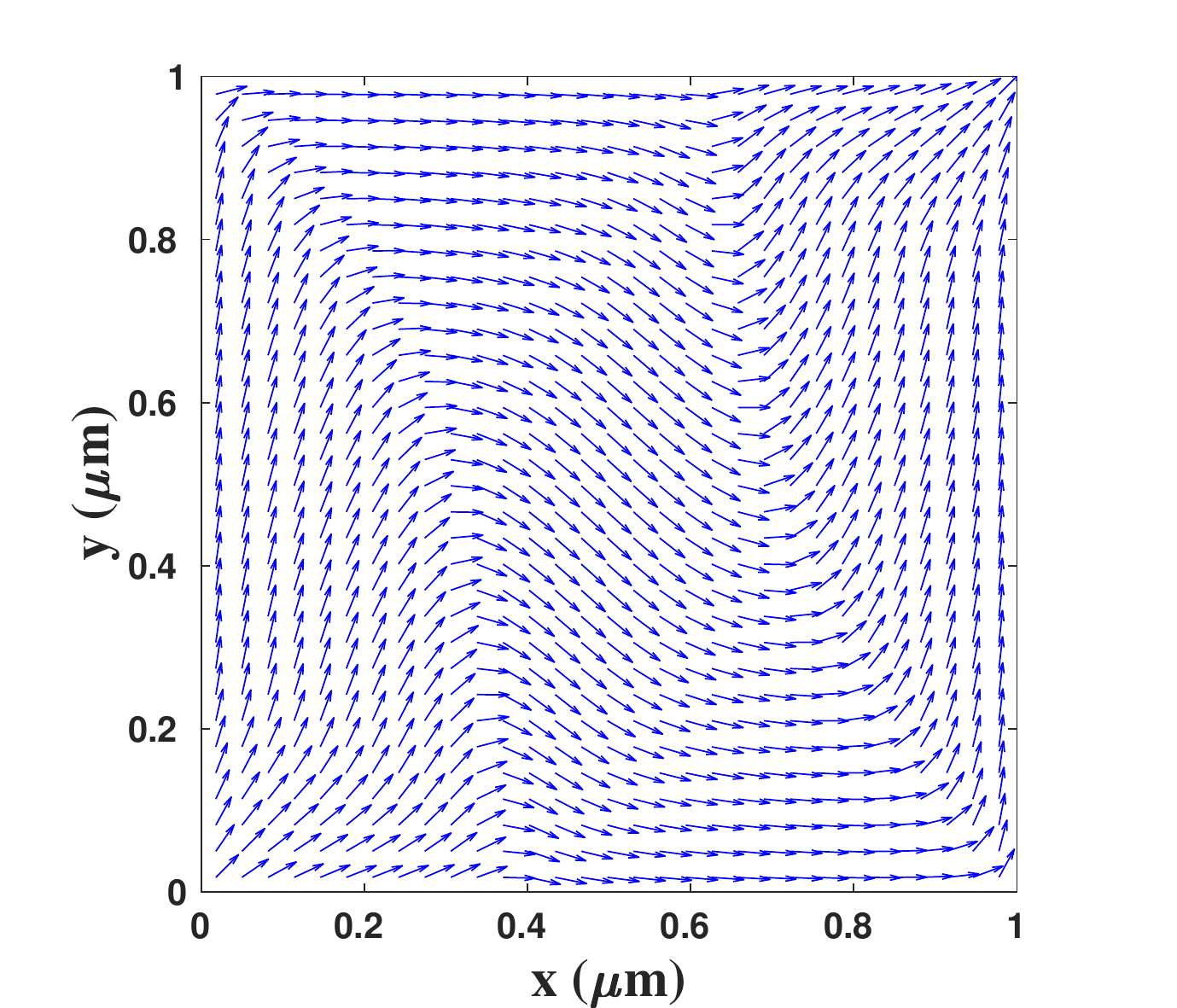}}
	\caption{Simulation of the full Landau-Lifshitz-Gilbert equation using Scheme B without any external field. The magnetization on the centered slice of the material in the $xy$ plane is used.
		Top row: $\alpha=0.1$; Bottom row: $\alpha=0.01$. Left column: a color plot of the angle between the in-plane magnetization and the $x$ axis;
		Right column: an arrow plot of the in-plane magnetization.}
	\label{IGSP}
\end{figure}

It is tested that GSPM in \cite{NumGSPM2001} was unstable with a very small damping parameter $\alpha$
and was resolved in \cite{ImproveGSPM2003}.
This section is devoted to the unconditional stability of Scheme A and Scheme B with respect to $\alpha$. We consider a thin film ferromagnet of size $1\;\mu\textrm{m}\times1\;\mu\textrm{m}\times0.02\;\mu\textrm{m}$
with the spatial gridsize  $4\;\textrm{nm}\times4\;\textrm{nm}\times 4\;\textrm{nm}$ and the temporal stepsize is $1$ picosecond.
Following \cite{ImproveGSPM2003}, we consider the full LLG equation with $\alpha=0.1$ and $\alpha=0.01$ and without the external field.
The initial state is $\m_0=(0,1,0)$ if $x\in[0,L_x/5]\cup[4L_x/5,L_x]$ and $\m_0=(1,0,0)$ otherwise. The final time is $10\;\textrm{ns}$.
In \Cref{GSPM,IGS,IGSP}, we present a color plot of the angle between the in-plane magnetization and the $x$ axis,
and an arrow plot of the in-plane magnetization for the original GSPM \cite{NumGSPM2001}, Scheme A, and Scheme B, respectively.
In these figures, $\alpha=0.1$ is presented in the top row and $\alpha=0.01$ is presented in the bottom row;
a color plot of the angle between the in-palne magnetization and the $x$ axis is presented in the left column
and  an arrow plot of the in-plane magnetization is presented in the right column.

\section{Conclusion}\label{sectionClusion}

In this paper, based on the original Gauss-Seidel projection methods,
we present two improved Gauss-Seidel projection methods with the first-order accuracy in time
and the second-order accuracy in space.
The first method updates the gyromagnetic term and the damping term simultaneously and follows
by a projection step, which requires to solve heat equations $5$ times at each time step.
The second method introduces two sets of approximate solutions, where
we update the gyromagnetic term and the damping term simultaneously for one set of approximate
solutions and apply the projection step to the other set of approximate solutions in an alternating manner.
Therefore, only $3$ heat equations are needed to be solved at each step.
Compared to the original Gauss-Seidel projection method, which solves heat equations $7$ times at each step,
savings of these two improved methods are about $2/7$ and $4/7$, which is verified by
both 1D and 3D examples for the same accuracy requirement.
In addition, unconditional stability with respect to both the grid size and the damping parameter is confirmed
numerically. Application of both methods to a realistic material is also presented with hysteresis loops
and magnetization profiles.

\section*{Acknowledgments}
This work is supported in part by the grants NSFC 21602149 (J.~Chen), NSFC 11501399 (R.~Du), the Hong Kong Research Grants Council (GRF grants 16302715,
16324416, 16303318 and NSFC-RGC joint research grant N-HKUST620/15) (X.-P.~Wang), and the Innovation Program for postgraduates in Jiangsu province via grant KYCX19\_1947 (C.~Xie).

\bibliographystyle{model1-num-names}
\bibliography{refs}

\end{document}